\documentclass[a4paper]{elsarticle}

\usepackage[utf8]{inputenc}
\usepackage{import}

\usepackage{fontawesome}
\usepackage{lipsum}

\usepackage[english]{babel}

\usepackage[margin=2.45cm]{geometry}

\usepackage{xcolor}

\usepackage{graphicx}
\usepackage{svg}
\usepackage{epstopdf}
\usepackage{tikz}
\usetikzlibrary{calc}
\usepackage{pgfplots}
\pgfplotsset{compat=newest}
\usepackage[percent]{overpic}
\usepackage{placeins}

\usepackage{circledsteps}
\pgfkeys{/csteps/fill color=white}
\pgfkeys{/csteps/inner ysep=5pt}
\pgfkeys{/csteps/inner xsep=5pt}

\usepackage{enumitem}
\setlist[1]{itemsep=-5pt}

\usepackage{amsmath}
\usepackage{amssymb,amsthm}
\usepackage{cases}
\usepackage{mathrsfs}
\usepackage{euscript}
\usepackage{bm}
\usepackage{stmaryrd}
\usepackage{mathtools}
\usepackage{empheq}

\theoremstyle{definition}

\usepackage{algorithm}
\usepackage{algpseudocode}

\usepackage{appendix}

\usepackage{pdflscape}
\usepackage{setspace}
\usepackage{framed}
\usepackage{pdfpages}

\usepackage{tabularx}
\usepackage{booktabs}

\usepackage{multicol}
\usepackage{multirow}
\usepackage{xparse}
\usepackage{xspace}
\usepackage{makecell}
\usepackage{tablefootnote}

\usepackage{footnote}

\usepackage{caption}
\usepackage{subcaption}

\usepackage[notref,notcite,final]{showkeys}

\usepackage[noprefix]{nomencl}
\makenomenclature

\renewcommand*\nompreamble{\begin{multicols}{2}}
\renewcommand*\nompostamble{\end{multicols}}

\usepackage{etoolbox}
\usepackage{stackengine}
\setstackgap{S}{1.75pt}

\biboptions{sort&compress}
\usepackage[numbers]{natbib}
\usepackage[colorlinks=true]{hyperref}
\usepackage{cleveref}

\def \usetodonotes{} 

\ifdefined\usetodonotes

\definecolor{C0}{HTML}{1F77B4}
\definecolor{C1}{HTML}{FF7F0E}
\definecolor{C2}{HTML}{2ca02c}
\definecolor{C3}{HTML}{d62728}
\definecolor{C4}{HTML}{9467bd}
\definecolor{C5}{HTML}{8c564b}

\usepackage[textsize=scriptsize]{todonotes}

\NewDocumentCommand{\JU}{O{} m}{%
    \todo[fancyline,color=C0!10,linecolor=C0,bordercolor=C0,#1]{JU: #2}%
}
\NewDocumentCommand{\ND}{O{} m}{%
    \todo[fancyline,color=C1!10,linecolor=C1,bordercolor=C1,#1]{ND: #2}%
}
\NewDocumentCommand{\SP}{O{} m}{%
    \todo[fancyline,color=C2!10,linecolor=C2,bordercolor=C2,#1]{SP: #2}%
}

\newlength{\todostripwidth}
\usepackage{atbegshi}
\AtBeginShipout{%
    \AtBeginShipoutUpperLeft{%
        \begin{tikzpicture}[overlay,remember picture]
            \fill [gray!0] (current page.north east) rectangle
                ($ (current page.south east) - (\todostripwidth,0cm) $);
        \end{tikzpicture}
    }
}

\else

\NewDocumentCommand{\JU}{O{} m}{}
\NewDocumentCommand{\ND}{O{} m}{}
\NewDocumentCommand{\SP}{O{} m}{}

\fi

\newcommand{\md}{\kern0.1em---\kern0.1em\ignorespaces}
\newcommand{\nd}{\texorpdfstring{\kern0.5pt--\kern0.5pt}{--}}

\makeatletter
\newcommand{\multiline}[1]{%
  \begin{tabularx}{\dimexpr\linewidth-\ALG@thistlm}[t]{@{}X@{}}
    #1
  \end{tabularx}
}
\makeatother

\definecolor{myblue}{RGB}{0, 0, 204}

\allowdisplaybreaks
\begin{document}

\begin{frontmatter}

\author[add1]{Sebastián Pazmiño}
\ead{pazminos@umich.edu}
\author[add2]{Nitin P.~Daphalapurkar}
\ead{nitin@lanl.gov}
\author[add1]{Jacinto Ulloa\corref{cor1}}
\ead{julloa@umich.edu}

\cortext[cor1]{Corresponding author}
\address[add1]{Department of Mechanical Engineering, University of Michigan, Ann Arbor, MI 48109, USA}
\address[add2]{Theoretical Division, Los Alamos National Laboratory, Los Alamos, NM 87545, USA}

\title{
Modeling fluid\nd solid interaction in three-dimensional granular systems: \\ Two-way coupling of SPH and LS-DEM
}

\begin{abstract}
    This paper presents a two-way coupled fluid\nd solid interaction framework for granular systems based on the level-set discrete element method (LS-DEM) and smoothed particle hydrodynamics (SPH). The proposed method is mesh-free, fully Lagrangian, and equipped with a seamless description of arbitrary solid shapes. In particular, a level-set representation of grain geometry is used for both contact detection and fluid\nd solid coupling, allowing for realistic, highly irregular grains to be handled within a unified computational framework. Moreover, the method incorporates a convenient treatment of moving boundaries and admits a simple computational implementation. The capabilities of the proposed SPH\nd LS-DEM framework are demonstrated through several numerical benchmarks involving dam-break flow, flow through porous media, granular collapse, and cyclic loading of a saturated granular specimen. The study sets the stage for several extensions to multiphysics problems across disciplines.
\end{abstract}

\begin{keyword}
Granular media; Fluid\nd solid interaction; Level-set discrete element method; Smoothed particle hydrodynamics; Multiphysics systems
\end{keyword}

\end{frontmatter}

\section{Introduction}
    
    Granular materials interacting with fluids arise in many natural and engineering systems. In geotechnical problems, the presence of interstitial fluid plays a central role in the transport~\cite{FreyChurch2009,Houssais2015}, erosion~\cite{moffat2011spatial,HunterBowman2018}, failure~\cite{rice1979earthquake,mihalache2014mathematical}, and liquefaction~\cite{Ishihara1993,andrade2009predictive} of granular soils. Coupled fluid\nd solid mechanics also governs large-scale geophysical processes, such as the deformation, fracture, and long-term evolution of sea ice~\cite{Herman2017,moncada2025simulating}. Moreover, in industrial systems, fluid\nd grain interactions are intrinsic to fluidized beds~\cite{AndersonJackson1967,Deen2007}, particle processing~\cite{Iveson2001,Claydon2022}, and clogging phenomena~\cite{MarinSouzy2025,Zhou2024}.

    Understanding and predicting these complex systems requires models that capture grain-scale contact dynamics and coupled fluid flow. Particle-resolved formulations approach this problem by solving the Navier\nd Stokes equations for the fluid phase while explicitly tracking the moving fluid\nd solid interfaces. Different strategies have been developed in this context, including immersed‑boundary methods with finite‑difference~\cite{uhlmann2005immersed} or finite‑volume~\cite{Breugem2012} schemes. Other examples are finite‑element formulations with fictitious domains~\cite{glowinski1999distributed} or embedded boundaries~\cite{Samaniego2015}. 
    
    For granular systems, a common approach is to represent the solid phase as rigid particles using the discrete element method (DEM), while accounting for fluid effects at varying levels of fidelity. The simplest DEM-based methods account for fluid effects without solving the Navier\nd Stokes equations for the fluid phase. Examples include the constant-volume DEM approach, where a global isochoric constraint is imposed to approximate undrained mechanical tests under fully saturated conditions~\cite{yimsiri2010analysis, o2011particulate, Martin2020}, and the coupled fluid method, which incorporates a Darcy-scale representation~\cite{liu2015numerical,Salimi2025}. Higher-fidelity methods, by contrast, explicitly model the fluid dynamics and fluid\nd solid interactions. Within this class, unresolved methods solve spatially averaged or volume‑filtered fluid equations and determine the forces acting on each grain from empirical or semi‑empirical drag laws. Examples include finite-volume and finite-element CFD\nd DEM formulations~\cite{tsuji1993discrete,CapecelatroDesjardins2013,Malipeddi2023,Coppin2023}, as well as unresolved schemes based on smoothed particle hydrodynamics (SPH)~\cite{Sun2013,Robinson2014,sizkow2021sph}. Conversely, resolved methods explicitly describe the flow around each grain, at the price of increased computational cost. Examples include fully resolved CFD\nd DEM formulations~\cite{Kloss2012,hager2014parallel,ZhangTahmasebi2022,hu2023insight,Hu2025,lai2023signed}, as well as approaches based on the lattice Boltzmann method (LBM)~\cite{feng2010combined,Galindo2013} and SPH~\cite{potapov2001liquid,peng2021fully,liu2022sph}. Hybrid methods have also been developed, which combine unresolved and resolved descriptions in different regions according to the local ratio between the grain size and the spatial resolution of the fluid~\cite{Chen2026}.

    Concerning the description of the solid phase, most early CFD\nd DEM contributions employed the traditional framework set forth by~\citet{cundall1979discrete} for spherical particles. However, the crucial role of particle shape in the emergent behavior of granular materials has been well established in both dry~\cite{cho2006,jerves2016effects,adesina2024study,seo2024simulation,zhao2023role} and fluid-coupled~\cite{ZhangTahmasebi2022,liu2022sph,Hu2025,sizkow2021sph,lai2023signed,Zhou2024,peng2021fully} systems. Several approaches exist in the DEM literature for incorporating realistic particle shapes, including sphere clumps~\cite{garcia2009b}, polyhedral particles~\cite{cundall1988,alonso2002}, spherical harmonics~\cite{zhou2017generation}, the potential particle approach~\cite{houlsby2009}, the granular element method~\cite{lim2014granular}, particles represented by surface meshes~\cite{zhan2021surface}, the level-set discrete element method (LS-DEM)~\cite{Kawamoto2016Jun}, and related signed-distance frameworks~\cite{lai2022}. Accordingly, a number of CFD\nd DEM contributions have employed non-spherical particles coupled to different fluid descriptions~\cite{ZhangTahmasebi2022,liu2022sph,Coppin2023,Hu2025,lai2023signed,sizkow2021sph,peng2021fully,Chen2026}. 
    
    In the present contribution, we focus on LS-DEM, which provides a seamless transition from three-dimensional images of real grains to high-fidelity virtual avatars~\cite{Kawamoto2016Jun,kawamoto2018}. This method has been applied across a wide range of physical systems~\cite{buarque2021,bustamante2020generalized,wang2021,zhou2023,moncada2025simulating,harmon2023predicting,Zhou2024} and combined with simplified fluid coupling strategies for particular problems~\cite{Zhou2024,moncada2023,moncada2025simulating}. Nevertheless, a general LS‑DEM‑based framework for fully coupled fluid\nd solid interaction has not yet been developed. We therefore present an SPH\nd LS‑DEM framework that achieves this coupling. We specifically consider SPH to exploit the natural compatibility between the two Lagrangian particle-based discretizations. Moreover, the meshfree character of SPH enables a seamless treatment of dynamic solid boundaries and is particularly advantageous for the highly irregular particle shapes enabled by LS‑DEM, for which mesh‑based CFD solvers would require conforming grids or embedded‑boundary treatments. Furthermore, in the present formulation, the same level-set field used for contact detection supplies the quantities required for fluid\nd solid coupling, greatly simplifying the computational~implementation. Several examples are presented to demonstrate the ability of the proposed framework to describe a wide range of fluid\nd solid interaction phenomena.

    The paper is structured as follows. Section~\ref{sec:govern} presents the mathematical formulation of LS-DEM and SPH. The coupling scheme and its computational implementation are then described in section~\ref{sec:coupling}. Finally, section~\ref{sec:numerical} presents the computational benchmarks, followed by conclusions and perspectives in section~\ref{sec:conc}.

\section{Governing equations and discretization}\label{sec:govern}

    This section introduces the governing equations describing the mechanical behavior of the solid and fluid phases, as well as the adopted numerical discretization methods and solution strategies. The solid phase consists of grains modeled as arbitrarily shaped rigid bodies governed by Newton\nd Euler dynamics, whose interactions are handled numerically using the level-set discrete element method (LS-DEM). On the other hand, the fluid phase is formulated at the continuum level and discretized via smoothed-particle hydrodynamics (SPH), which is well-suited for time-dependent boundaries around irregular moving solids. The formulation explicitly includes fluid-induced forces and moments on the grains and the reciprocal effect of moving grains on the fluid, with the interaction terms defined through the coupling scheme in section~\ref{sec:coupling}.

    \subsection{Physical setting and domain decomposition}

    Consider a granular assembly composed of $N^{\mathrm g}$ solid grains interacting with a fluid, evolving over a time interval $\mathcal{T}=[0,t_\mathrm{max}]$. At time $t\in\mathcal{T}$, the system occupies a spatial domain $\Omega(t)\subset\mathbb{R}^n$, with $n=3$. Letting $\Omega_i^{\mathrm g}(t)\subset\mathbb{R}^n$ denote the spatial domain occupied by grain $i$, the granular phase follows as the union of all grain domains,
    \begin{equation}
        \Omega^{\mathrm g}(t) \coloneqq \bigcup_{i=1}^{N^{\mathrm g}} \Omega_i^{\mathrm g}(t).
    \end{equation}
    The fluid domain is the complement of the granular phase in $\Omega(t)$,
    \begin{equation}
        \Omega^{\mathrm f}(t) \coloneqq \Omega(t) \setminus \Omega^{\mathrm g}(t),
    \end{equation}
    so that $\Omega(t)=\Omega^{\mathrm g}(t)\cup\Omega^{\mathrm f}(t)$ and $\Omega^{\mathrm g}(t)\cap\Omega^{\mathrm f}(t)=\varnothing$. Defining $\Gamma^{\mathrm g}(t)=\partial\Omega^{\mathrm g}(t)$ and $\Gamma^{\mathrm f}(t)=\partial\Omega^{\mathrm f}(t)$, the fluid\nd grain interface follows as the common boundary,
    \begin{equation}\label{eq:fsint}
        \Gamma^{\mathrm{fg}}(t) \coloneqq \Gamma^{\mathrm g}(t)\cap\Gamma^{\mathrm f}(t).
    \end{equation}

    \subsection{Solid phase}\label{sec:solid}

        \subsubsection{Kinematics and mechanical balance}
            
            The solid phase consists of $N^{\mathrm g}$ arbitrarily shaped rigid grains whose motion is described by translational and rotational degrees of freedom for all $t\in \mathcal{T}$. Specifically, for each grain $i\in\{1,\dots,N^{\mathrm g}\}$, let $\bm x_i(t)\in\mathbb{R}^n$ and $\bm\omega_i(t)\in\mathbb{R}^{n(n-1)/2}$ denote the centroid position and angular velocity, respectively. The grain orientation is represented by a unit quaternion $\bm q_i(t)\in\mathbb{R}^4$, encoding the rotation from the body-fixed frame of grain $i$ to the global frame. We denote the corresponding rotation tensor by $\mathbf R_i(t)\in\mathrm{SO}(n)$. 
            
            Assuming rigid-body kinematics, the dynamics of each grain are governed by the Newton\nd Euler equations of motion, i.e., for all $i\in\{1,\dots,N^\mathrm{g}\}$ and $t\in\mathcal{T}$,
            \begin{subequations}\label{eq:solid_dynamics}
                \begin{align}
                m_i\ddot{\bm x}_i + \xi m_i\dot{\bm x}_i &= \sum_{\alpha\in\mathcal{C}_i}{\bm f}^\alpha_i + \bm{f}_i^\mathrm{ext} + \bm{f}_i^\mathrm{fluid},\label{eq:solid_trans}\\[0.5em]
                \mathbf{J}_i \dot{\bm{\omega}}_i
                +
                \bm{\omega}_i \times (\mathbf{J}_i \bm{\omega}_i)
                +
                \xi \, \mathbf{J}_i \bm{\omega}_i
                &=
                \sum_{\alpha \in \mathcal{C}_i} \bm{m}_i^\alpha
                +
                \bm{m}_i^{\mathrm{ext}}
                +
                \bm{m}_i^{\mathrm{fluid}}, \label{eq:solid_rot}
                \end{align}
            \end{subequations}
            with the kinematic relation
            \begin{equation}\label{eq:quatkin}
                \dot{\bm q}_i(t)=\frac12\mathbf Q(\bm q_i(t))\,\bm\omega_i(t),
            \end{equation}
            where $\mathbf Q(\bm q_i)\in\mathbb{R}^{4\times 3}$ is the standard quaternion kinematic matrix,
            and initial conditions 
            \begin{equation}\label{eq:IC}
                \bm x_i(0)=\bm x_{i,0}, \qquad \dot{\bm x}_i(0)=\bm v_{i,0}, \qquad
                \bm q_i(0)=\bm q_{i,0}, \qquad {\bm\omega}_i(0)=\bm\omega_{i,0}.
            \end{equation}
            In equations~\eqref{eq:solid_dynamics}, $m_i$ is the mass of grain $i$ and $\mathbf J_i$ is the inertia tensor about the centroid expressed in the global frame.  On the right-hand side, the first terms represent the resultant of all forces ${\bm f}^\alpha_i$ and moments ${\bm m}^\alpha_i$, $\alpha\in\mathcal{C}_i$, acting on the grain centroid due to contacts with other grains or walls, with $\mathcal{C}_i$ denoting a set of contact indices. The second terms represent long-range, prescribed external forces $\bm{f}_i^\mathrm{ext}$ (e.g., for gravity, $\bm{f}_i^\mathrm{ext}= m_i\,\bm g$) and moments $\bm{m}_i^\mathrm{ext}$ (usually zero, unless additional physics are present). Finally, the third terms $\bm{f}_i^{\mathrm{fluid}}$ and $\bm{m}_i^{\mathrm{fluid}}$ represent fluid-induced forces and moments, to be specified in section~\ref{sec:coupling}. The parameter $\xi\ge 0$ is a numerical damping coefficient introduced for stabilization.

            Once the forces and moments on the right-hand side of~\eqref{eq:solid_dynamics} are known, the translational and rotational states $(\bm x_i,\bm q_i)$ of each grain are updated from~\eqref{eq:solid_dynamics}\nd \eqref{eq:IC} using standard time integration procedures (see, e.g.,~\citet{lim2014granular} and~\citet{Kawamoto2016Jun} for details). It thus remains to resolve the granular interactions, namely the contact forces ${\bm f}^\alpha_i$ and moments ${\bm m}^\alpha_i$ in~\eqref{eq:solid_dynamics}, for which LS-DEM is employed.
            
        \subsubsection{Interaction forces: level-set discrete element method}
    
             The cornerstone of LS-DEM~\cite{Kawamoto2016Jun} is the use of level-set functions to characterize the morphology of each grain. Specifically, each grain $j$ occupies a fixed domain $\hat{\Omega}_j^{\mathrm g}\subset\mathbb R^n$ with boundary $\hat{\Gamma}_j^{\mathrm g}$, defined in a body-fixed coordinate system centered at the grain's centroid. The shape is then encoded by the time-independent signed-distance function
            \begin{equation}\label{eq:lset}
                \varphi_j(\bm\xi)
                =
                \begin{cases}
                    \displaystyle-\min_{\bm\eta\in\hat{\Gamma}_j^{\mathrm g}}\|\bm\xi-\bm\eta\|,
                    & \bm\xi\in\hat{\Omega}_j^{\mathrm g},\\[0.75em]
                    \displaystyle+\min_{\bm\eta\in\hat{\Gamma}_j^{\mathrm g}}\|\bm\xi-\bm\eta\|,
                    & \bm\xi\notin\hat{\Omega}_j^{\mathrm g},
                \end{cases}
            \end{equation}
            where $\bm\xi$ denotes a body-fixed coordinate and $\hat{\Gamma}_j^{\mathrm g}=\{\bm\xi \;:\; \varphi_j(\bm\xi)=0\}$. This representation accommodates arbitrary grain shapes, enabling a direct connection with experimental images~\cite{Kawamoto2016Jun}. The current grain domain and boundary follow from the centroid $\bm x_j(t)$ and rotation tensor $\mathbf R_j(t)\in\mathrm{SO}(n)$,
            \begin{equation}
                \Omega_j^{\mathrm g}(t)
                =
                \bm x_j(t)+\mathbf R_j(t)\hat{\Omega}_j^{\mathrm g},
                \qquad
                \Gamma_j^{\mathrm g}(t)
                =
                \bm x_j(t)+\mathbf R_j(t)\hat{\Gamma}_j^{\mathrm g}.
            \end{equation}
            Accordingly, a spatial point $\bm x$ is mapped to the body-fixed frame of grain $j$ as
            \begin{equation}
                \bm\xi_j(\bm x,t)
                =
                \mathbf R_j^{\mathsf T}(t)
                \bigl(\bm x-\bm x_j(t)\bigr).
            \end{equation}
            Thus, level-set evaluations are performed as $\varphi_j(\bm\xi_j(\bm x,t))$, while the stored level-set function remains fixed. The explicit time dependence of the rigid-body quantities is omitted below for conciseness. 
            
            Contact between a grain and another object\md such as a neighboring grain or a wall\md is determined by mapping surface points of the other object to the grain body-fixed frame and evaluating the signed-distance function~\eqref{eq:lset}. In practice, the signed-distance function for each grain is discretized over a grid, with evaluations performed using bi-linear interpolation in 2D or tri-linear interpolation in 3D.

            \begin{figure}[t!]
                \centering
                \vspace{1em}
                  \includegraphics[width=0.8\textwidth]{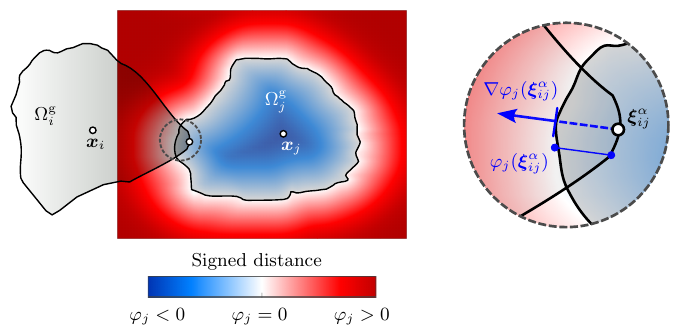}
                  \caption{Two-dimensional illustration of LS-DEM, showing contact between grains $i$ and $j$ using the level-set function $\varphi_j$.}
                \vspace{1em}
                  \label{fig:lsdem}
            \end{figure}

            Suppose that contact is detected between grain $j$ and a point $\bm x_i^\alpha\in\Gamma_i^{\mathrm g}$ on the surface of grain $i$ (Figure~\ref{fig:lsdem}). The point is mapped to the body frame of grain $j$ as
            \begin{equation}
                \bm\xi_{ij}^\alpha
                =
                \mathbf R_j^{\mathsf T}
                \bigl(\bm x_i^\alpha-\bm x_j\bigr).
            \end{equation}
            The penetration distance and spatial outward normal associated with grain $j$ are then
            \begin{equation}
                \delta_i^\alpha
                =
                \max\{-\varphi_j(\bm\xi_{ij}^\alpha),0\},
                \qquad
                \bm n_i^\alpha
                =
                \mathbf R_j
                \frac{\nabla\varphi_j(\bm\xi_{ij}^\alpha)}
                {\|\nabla\varphi_j(\bm\xi_{ij}^\alpha)\|}.
            \end{equation}
            If $\delta_i^\alpha>0$, the contact force at $\bm x_i^\alpha$ is decomposed as $\bm{f}^\alpha_i = \bm{f}^{\alpha,\mathrm{n}}_{i} + \bm{f}^{\alpha,\mathrm{s}}_{i}$, where, assuming linear elastic contacts with Coulomb friction, the normal and shear components read
            \begin{equation}
            \label{eq:contactForce}
            \begin{aligned}
            & \quad \bm{f}^{\alpha,\mathrm{n}}_{i} = k^{\vphantom{c}}_\mathrm{n}\delta_i^\alpha\bm{n}^\alpha_i, \qquad
            &\begin{cases}\bm{f}^{\alpha,\mathrm{s}}_{i} \leftarrow \mathbf{R}^\alpha_i\bm{f}^{\alpha,\mathrm{s}}_{i}  - k_\mathrm{s}\,\Delta\bm{s}^\alpha_i, \\[0.75em] 
            \bm{f}^{\alpha,\mathrm{s}}_{i} \leftarrow  \dfrac{\bm{f}^{\alpha,\mathrm{s}}_{i}}{\Vert\bm{f}^{\alpha,\mathrm{s}}_{i}\Vert} \min\big\{ \Vert\bm{f}^{\alpha,\mathrm{s}}_{i}\Vert, \; \mu \, \Vert \bm{f}^{\alpha,\mathrm{n}}_{i}\Vert \big\}, \end{cases}
            \end{aligned}
            \end{equation}
            where $\Delta\bm{s}^\alpha_i$ is the incremental tangential displacement and $\mathbf{R}^\alpha_i$ rotates the previous contact normal to the current direction. The constants $k_\mathrm{n}$, $k_\mathrm{s}$, and $\mu$ denote the normal stiffness, shear stiffness, and friction coefficient, respectively. Moments acting on the grain centroid are updated similarly, 
            \begin{equation}\label{eq:momentni}
            \bm{m}^\alpha_i=\bm{r}_i^\alpha\times \bm{f}^\alpha_i,
            \end{equation}        
            where $\bm{r}_i^\alpha=\bm{x}_i^\alpha-\bm{x}_i$ is the branch vector connecting the centroid of grain $i$ to the surface contact point. We refer the reader to~\citet{Kawamoto2016Jun} for further details on LS-DEM and its implementation.
                        
    \subsection{Fluid phase}\label{sec:fluid}

        \subsubsection{Governing equations}
            The fluid is modeled as a continuum occupying $\Omega^{\mathrm f}(t)\subset\mathbb{R}^n$, with boundary $\Gamma^{\mathrm f}(t)$. The Dirichlet and Neumann parts of the external fluid boundary are denoted by $\Gamma_{\mathrm D}^{\mathrm f}(t)$ and $\Gamma_{\mathrm N}^{\mathrm f}(t)$, respectively. Together with the fluid\nd solid interface $\Gamma^{\mathrm{fg}}(t)$ defined in~\eqref{eq:fsint}, the fluid boundary is decomposed as
            \begin{equation}
                \Gamma^{\mathrm f}(t)
                =
                \Gamma_{\mathrm D}^{\mathrm f}(t)
                \cup
                \Gamma_{\mathrm N}^{\mathrm f}(t)
                \cup
                \Gamma^{\mathrm{fg}}(t).
            \end{equation}
    
            Fluid evolution is characterized by the velocity field $\bm v(\bm x,t)\in\mathbb{R}^n$ and the pressure field $p(\bm x,t)\in\mathbb{R}$. We assume incompressible flow with constant density $\rho>0$. The behavior of the fluid is then governed by the incompressibility constraint and the balance of linear momentum; that is, for all $t\in\mathcal{T}$,
            \begin{subequations}\label{eq:inc_fleq}
                \begin{align}
                &\mathrm{div}\,\bm v=0
                &&\hspace{-8em}\text{in } \Omega^{\mathrm f}(t), \label{eq:inc_cont}\\
                &\rho\frac{\mathrm D\bm v}{\mathrm Dt}
                =\mathrm{div}\,\bm\sigma+\rho\bm b
                &&\hspace{-8em}\text{in } \Omega^{\mathrm f}(t), \label{eq:inc_mom}
                \end{align}
            \end{subequations}
            where ${\mathrm D\bm v}/{\mathrm Dt}=\partial\bm v/\partial t +(\nabla\bm v)\bm v$ is the material acceleration in Eulerian form, $\bm\sigma$ is the Cauchy stress tensor, and $\bm b$ is a prescribed body force per unit mass. Assuming a Newtonian constitutive response, the stress tensor is given by
            \begin{equation}
                \bm\sigma=-p\,\mathbf{I}+2\eta\,\mathrm{sym}\,\nabla\bm v,
                \label{eq:ns_stress}
            \end{equation}
            where $\eta>0$ is the dynamic viscosity and $\mathbf{I}$ is the second-order identity tensor. The system is supplemented with boundary conditions for all $t\in\mathcal T$,
            \begin{subequations}
                \begin{align}
                &\bm v=\bar{\bm v}
                &&\hspace{-8em}\text{on } \Gamma_{\mathrm D}^{\mathrm f}(t), \label{eq:inc_bcD}\\
                &\bm\sigma\bm n=\bar{\bm t}
                &&\hspace{-8em}\text{on } \Gamma_{\mathrm N}^{\mathrm f}(t), \label{eq:inc_bcN}
                \end{align}
            \end{subequations}
            where $\bm n$ denotes the outward unit normal to $\Gamma^{\mathrm f}(t)$, $\bar{\bm v}$ is a prescribed velocity on $\Gamma_{\mathrm D}^{\mathrm f}(t)$, and $\bar{\bm t}$ is a prescribed traction on $\Gamma_{\mathrm N}^{\mathrm f}(t)$. Conditions along $\Gamma^{\mathrm{fg}}(t)$ are resolved through the coupling scheme introduced in section~\ref{sec:coupling}. Finally, the system is supplemented with initial conditions,
            \begin{equation}
                \bm v(\cdot,0)=\bm v_0
                \qquad\text{in } \Omega^{\mathrm f}(0).
            \end{equation}
      
        \subsubsection{Smoothed-particle hydrodynamics}

            Equations~\eqref{eq:inc_cont}\nd \eqref{eq:inc_mom} are discretized using a Lagrangian SPH formulation, which is convenient for handling the moving fluid\nd solid interfaces induced by grain motion. In SPH, the fluid is represented by a finite set of $N^{\mathrm{f}}$ moving particles carrying mass and state variables, e.g., $\rho_a(t)$, $p_a(t)$, or $\bm v_a(t)$ for all $a\in\{1,\dots,N^{\mathrm f}\}$, where $(\Box)_a$ denotes the evaluation of a spatial field $\Box$ at a discrete particle location $\bm x_a(t)$. The governing equations~\eqref{eq:inc_fleq}\nd \eqref{eq:ns_stress} are then enforced discretely at $\bm x_a(t)$. To this end, SPH employs kernel-based approximations of the spatial fields and their differential operators~\cite{Monaghan2012}. 
            
            Let us briefly recall this concept for the reader's convenience. Consider a sufficiently smooth scalar field $f:\mathbb{R}^n\to\mathbb{R}$. We may express an evaluation $f(\bm x)$ using a regularization of the Dirac-$\delta$ identity,
            \begin{equation}\label{eq:kernel_regularization}
                f(\bm x)
                =
                \int_{\mathbb R^n}f(\bm y)\,\delta(\bm x-\bm y)\,\mathrm d\bm y
                \approx
                \int_{\mathbb R^n}f(\bm y)\,W(\bm x-\bm y,h)\,\mathrm d\bm y,
            \end{equation}
            where $W(\cdot,\,h):\mathbb{R}^n\to\mathbb{R}$ is a smoothing kernel with characteristic length $h$. The kernel is assumed to satisfy the standard normalization and consistency conditions
            \begin{equation}\label{eq:kernelAprox}
                \int_{\mathbb{R}^n} W(\bm r,h)\,\mathrm d\bm r = 1,
                \qquad
                \lim_{h\to 0} W(\bm r,h)=\delta(\bm r).
            \end{equation}
            Introducing the differential mass $\mathrm d m(\bm y)=\rho(\bm y)\,\mathrm d\bm y$ yields the continuous (integral) SPH form
            \begin{equation}\label{eq:continuousForm}
                f(\bm x)
                \approx
                \int_{\mathbb R^n}
                \frac{f(\bm y)}{\rho(\bm y)}
                W(\bm x-\bm y,h)\,
                \mathrm dm(\bm y),
            \end{equation}
            which, upon particle discretization, leads to the standard discrete approximation
            \begin{equation}\label{eq:discretization}
                f_a  = f(\bm x_a) \approx \sum_{b\in \mathcal H_a} \frac{m_b}{\rho_b}\,f(\bm x_b)\,W(\bm x_a-\bm x_b,h),
            \end{equation}
            where $\mathcal H_a = \mathcal H(\bm x_a)=\{b\in\{1,\dots,N^{\mathrm f}\} \; : \; W(\bm x_a-\bm x_b,h)\neq 0\}$ is an index set of fluid particles within the compact support of the kernel at $\bm x_a$ (set by the smoothing length $h$). Accordingly, the differential operators $\nabla$ and $\nabla^2$ take the discrete forms 
            \begin{align}\label{eq:gradientDiscrete}
                (\nabla f)_a = \nabla f(\bm x_a) \approx \sum_{b\in \mathcal H_a} \frac{m_b}{\rho_b}\,f(\bm x_b)\,\nabla W(\bm x_a-\bm x_b,h),\\
            \label{eq:laplacianDiscrete}
                (\nabla^{2} f)_a = \nabla^{2} f(\bm x_a) \approx \sum_{b\in \mathcal H_a} \frac{m_b}{\rho_b}\,f(\bm x_b)\,\nabla^{2} W(\bm x_a-\bm x_b,h).
            \end{align}

\begin{figure}[t!]
            \vspace*{1em}
            \begin{center}
            \begin{tabular}{cc}
            \begin{overpic}[width=0.46\textwidth]{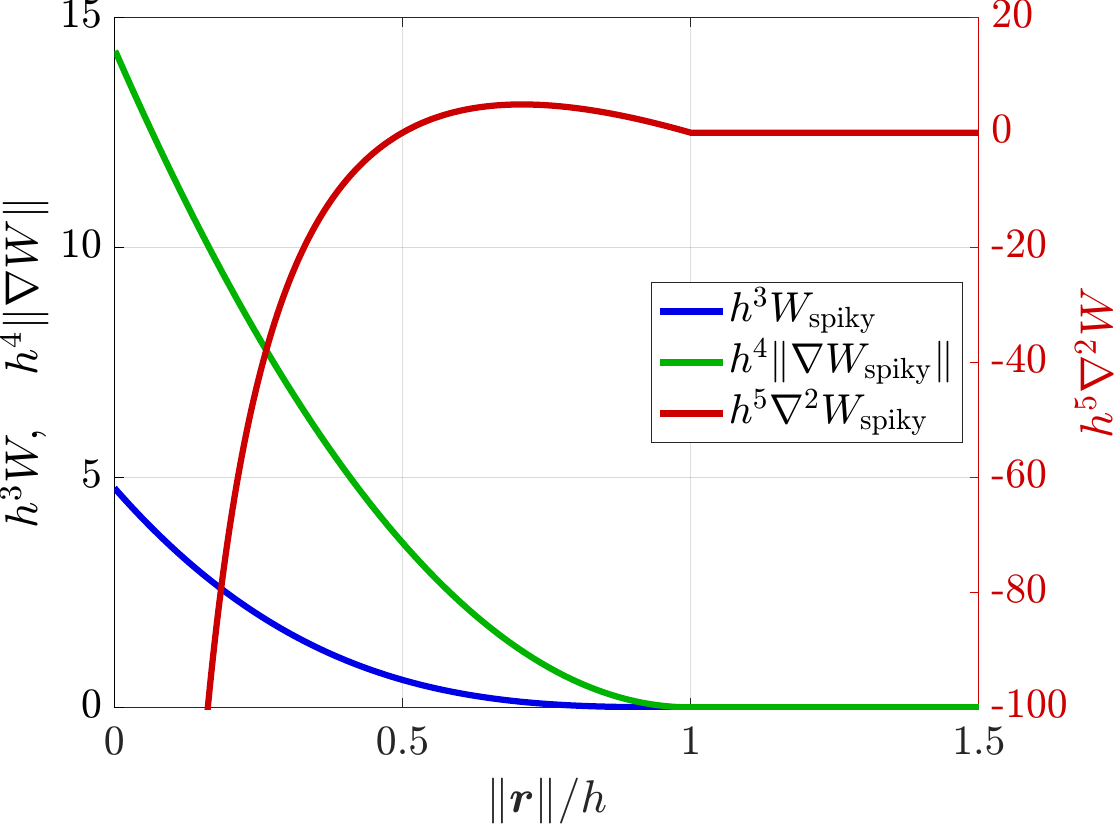}
                \put(0,2){\small{(a)}}
            \end{overpic} &
            \begin{overpic}[width=0.46\textwidth]{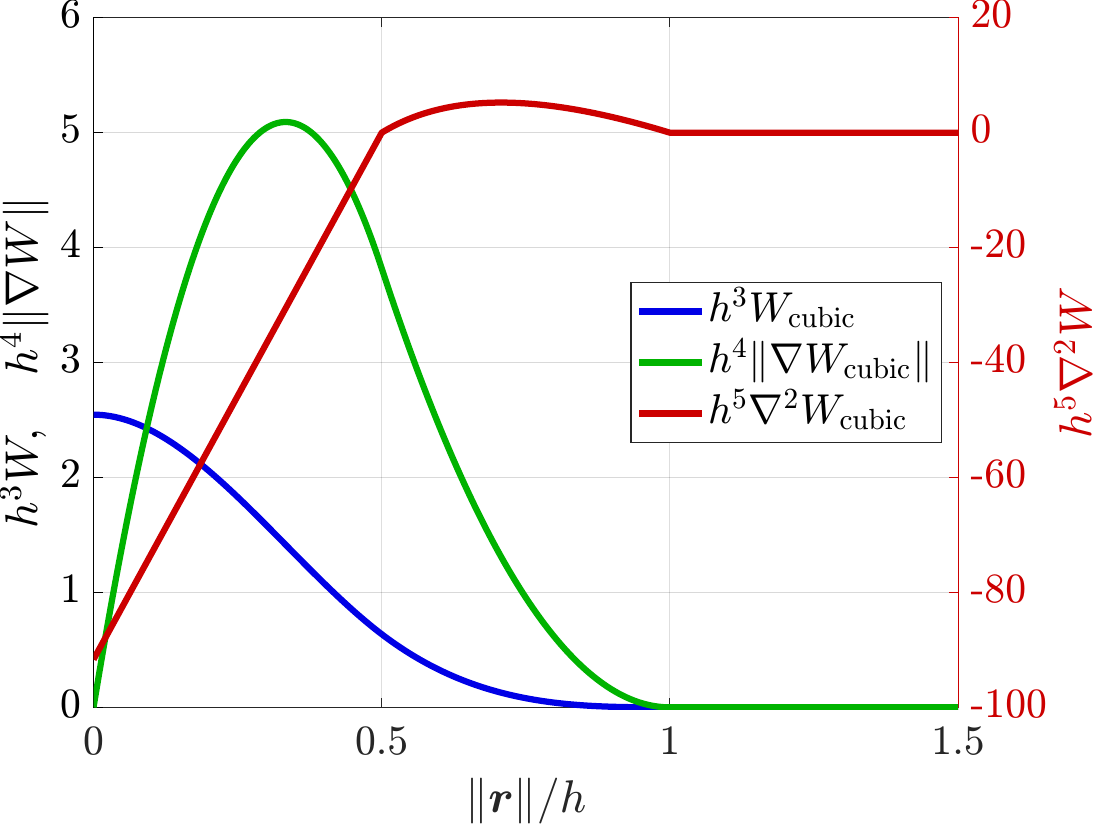}
                \put(0,2){\small{(b)}}
            \end{overpic}
            \end{tabular}
            \end{center}
            \vspace*{-1em}
            \caption{Kernel functions and derivatives for (a) the spiky kernel and (b) the cubic kernel normalized by $h$.}
        \label{fig:kernels}
        \end{figure}

        It remains to define the kernel function, whose choice affects both accuracy and computational cost. For instance, a Gaussian kernel has infinite support, leading to an interaction cost of $\mathcal O\!\left((N^{\mathrm f})^2\right)$. However, in most hydrodynamic processes, long-range contributions are negligible. Thus, compact-support kernels are typically employed, reducing the cost to $\mathcal O(N^{\mathrm f}N_W)$, where $N_W$ denotes the typical number of neighbors within the kernel support~\cite{sigalotti2021sph}. An example is the spiky kernel, originally proposed by~\citet{desbrun1996smoothed} and later popularized for fluid simulations by~\citet{Muller2003}. Its non-vanishing gradient near the origin is beneficial for pressure-gradient computations and helps mitigate particle overlap. The spiky kernel is defined in three dimensions with normalization constant $c_{n}=15/(\pi h^6)$ as
        \begin{equation}
        \label{eq:spikyK}
            W_{\mathrm{spiky}}(\bm r,h)=
            \begin{dcases}
            c_{n}\,(h-\|\bm r\|)^3 &\text{if} \ \  0\le \|\bm r\|\le h,\\[0.3em]
            0 & \text{otherwise}.
        \end{dcases}
        \end{equation}
        Alternatively, a typical choice in the SPH literature is the cubic spline kernel~\cite{Koschier2020arXiv}, with continuous first and second derivatives providing numerical stability and accuracy for general field interpolations, especially when calculating density or viscous forces. Letting $q \coloneqq \Vert\bm{r}\Vert/h$, the cubic kernel is defined in three dimensions with normalization constant $c_{n} = 8/(\pi h^{3})$ as
        \begin{equation}
            \label{eq:cubic}
            W_{\rm cubic}(\bm{r},h) = c_n 
            \begin{cases}
            6\,(q^{3} -q^{2}) + 1 &\text{if} \ \ 0 \le q \le \tfrac{1}{2},\\[0.3em]
            2\,(1-q)^{3}  &\text{if} \ \ \tfrac{1}{2} < q \le 1,
            \\[0.3em]
            0 & \text{otherwise}.
            \end{cases}
        \end{equation}
        Figure~\ref{fig:kernels} shows both kernels and their derivatives. In this work, we adopt the spiky kernel~\eqref{eq:spikyK} to assess the stability of our model under free-flow conditions and otherwise employ the cubic kernel~\eqref{eq:cubic}.

        \subsubsection{Discrete fluid equations}

            Let us now apply the preceding discretization method to the fluid equations~\eqref{eq:inc_cont}\nd \eqref{eq:inc_mom}. In SPH, incompressibility may be enforced directly through the constraint~\eqref{eq:inc_cont} or approximated within a weakly-compressible SPH (WCSPH) setting. In the latter, incompressibility is enforced approximately by allowing small density variations and prescribing the pressure explicitly as a function of density, $p_a=p(\rho_a)$, such that $\rho_a$ remains close to a reference value $\rho_0$. In the present work, we adopt WCSPH and thus evaluate the particle density from~\eqref{eq:discretization}, i.e., for each $a\in\{1,\dots,N^{\mathrm f}\}$,
            \begin{equation}\label{eq:discreteDensitya}
                \rho_a \approx \sum_{b\in\mathcal H_a} m_b\,W(\bm x_a-\bm x_b,h).
            \end{equation}
            The corresponding pressure law is taken as~\cite{WANG20169625}
            \begin{equation}\label{eq:pressurea}
                p_a = k_\mathrm{e} \max(\rho_a-\rho_0,0),
            \end{equation}
            where $k_\mathrm{e}$ is a stiffness parameter controlling density fluctuations. The truncation $p_a\ge 0$ is used to mitigate non-physical particle cohesion under negative pressure~\cite{monaghan2000sph}.

            With $\rho_a$ and $p_a$ defined, we next discretize the balance of linear momentum~\eqref{eq:inc_mom} at the particle locations, using the constitutive relation~\eqref{eq:ns_stress}. We then have, for all $a\in\{1,\dots,N^{\mathrm f}\}$,
            \begin{equation}
                \label{eq:momentumFluid}            
                \bm a_a = \frac{\mathrm D \bm v_a}{\mathrm D t}
                =
                -\left(\frac{1}{\rho}\nabla p\right)_{\!a}
                \,+\, \nu\left(\nabla^{2}\bm v\right)_{\!a}
                \,+\, \bm b_a,
            \end{equation}
            where $\nu=\eta/\rho_0$.             Equation~\eqref{eq:momentumFluid} governs the velocity evolution along the SPH particle trajectories. Once the acceleration contributions are specified, the particle states $(\bm x_a,\bm v_a)$ can be advanced in time using standard explicit integration schemes. The discrete approximations of the acceleration contributions on the right-hand side of~\eqref{eq:momentumFluid} are described below.
            
            In equation~\eqref{eq:momentumFluid}, the pressure-gradient $\nabla p$ is not discretized directly from~\eqref{eq:gradientDiscrete}, as this may lead to non-conservative pairwise interactions. Instead, following references~\cite{monaghan1992sph,price2012smoothed}, we take
            \begin{equation}\label{eq:pressureAccelerationdiscrete}
                \left(\frac{1}{\rho}\nabla p\right)_{\!a}
                \;\approx\;
                \sum_{b\in\mathcal H_a} m_b
                \left(
                \frac{p_a}{\rho_a^{2}}+\frac{p_b}{\rho_b^{2}}
                \right)\nabla W(\bm x_a-\bm x_b,h).
            \end{equation}
            The viscous contribution in~\eqref{eq:momentumFluid} involves $\nabla^2\bm v$, which may be directly evaluated from~\eqref{eq:laplacianDiscrete}. However, this straightforward discretization is sensitive to particle disorder and may be numerically unstable as $h$ approaches zero. We therefore employ the Brookshaw-type discretization~\cite{monaghan2005sph}, i.e., taking the first gradient via~\eqref{eq:gradientDiscrete} and the second via finite differences, such that
            \begin{equation}\label{eq:viscosityAcceleration}
                \nu\,(\nabla^{2}\bm v)_a
                \;\approx\;
                2\nu\,(n+2)\sum_{b\in\mathcal H_a}\frac{m_b}{\rho_b}
                \left(
                \frac{(\bm v_a-\bm v_b)\cdot(\bm x_a-\bm x_b)}
                {\|\bm x_a-\bm x_b\|^2+\epsilon h^2}
                \right)\nabla  W(\bm x_a-\bm x_b,h),
            \end{equation}
            where $\epsilon h^2$ prevents division by zero, with typical values $\epsilon\sim 0.01$~\cite{price2012smoothed,monaghan2005sph}.

            \begin{figure}[t!]
                \centering
                \includegraphics[width=0.7\textwidth]{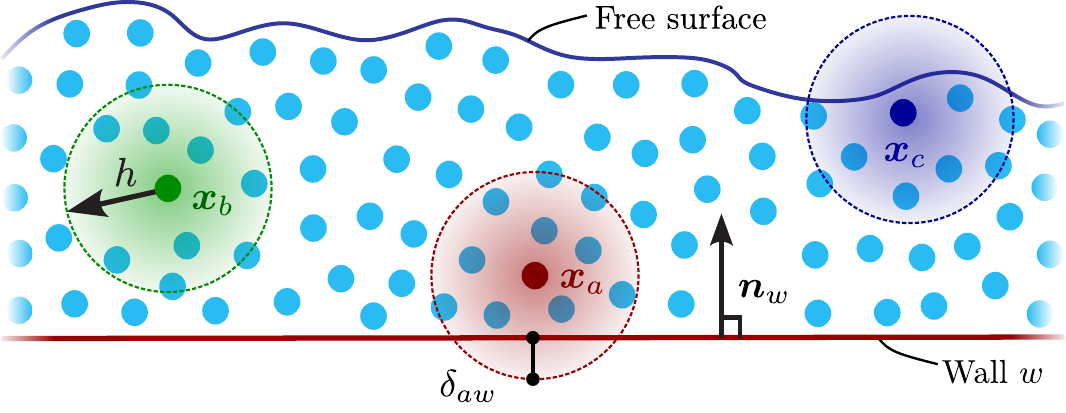}
                \caption{Schematic representation of the SPH kernel support with radius $h$ at an interior particle location $\bm x_b$ (green); a particle location $\bm x_a$ near a wall $w$, with normal $\bm n_w$ and overlap $\delta_{aw}$ (red); and a particle location $\bm x_c$ near the free surface, illustrating neighbor deficiency within the kernel support  (blue).}
                \label{fig:BoundaryConditionsWalls}
            \end{figure}      
     
            The remaining term in~\eqref{eq:momentumFluid} collects external forces and interaction effects, decomposed as
            \begin{equation}\label{eq:b_decomp}
                \bm b_a = \sum_{w=1}^{N^{\mathrm w}} \bm b^{\mathrm{wall}}_{aw} + \bm b_a^{\mathrm{ext}}  + \bm b_a^{\mathrm{free}} + \bm b_a^{\mathrm{solid}}.
            \end{equation}
             Here, $\bm b_a^{\mathrm{ext}}$ denotes prescribed body accelerations (e.g., for gravity, $\bm b_a^{\mathrm{ext}}=\bm g$). Boundary conditions are imposed through particle\nd wall interactions $\bm b_{aw}^{\mathrm{wall}}$, considering a container boundary composed of $N^{\mathrm w}$ planar walls with inward unit normals $\bm n_w$ and signed distances $d_{aw}$ from particle $a$ to wall $w$ (Figure~\ref{fig:BoundaryConditionsWalls}). This implicit wall representation avoids introducing additional boundary particles and provides an efficient distance-based contact evaluation. Assuming a linear elastic response, the acceleration due to wall forces reads
            \begin{equation}\label{eq:b_wall}
               \bm b^{\mathrm{wall}}_{aw}
                =\frac{1}{m_a} k_{{\mathrm n}w}\delta_{aw}\bm n_w, \qquad  \delta_{aw}=\max(-d_{aw},0),
            \end{equation}
            where $k_{\mathrm{n}w}$ is the particle\nd wall stiffness constant, set as $k_{\mathrm{n}w} = 40$~N/m throughout this work.

            The term $\bm b_a^{\mathrm{free}}$ accounts for surface-tension effects, introducing an inter-particle cohesive correction to mitigate the density underestimation in~\eqref{eq:discreteDensitya} and the resulting pressure underestimation in~\eqref{eq:pressurea} at free surfaces caused by neighbor deficiency  (Figure~\ref{fig:BoundaryConditionsWalls}). The corresponding acceleration reads~\cite{WANG20169625}
            \begin{equation}
                \bm b_a^{\mathrm{free}}=\frac{1}{m_a}\sum_{b\in\mathcal H_a\setminus\{a\}} \!\! s
                \cos\!\left(\frac{1.5\pi}{3h}\,\|\bm x_b-\bm x_a\|\right)
                \frac{\bm x_b-\bm x_a}{\|\bm x_b-\bm x_a\|},
            \end{equation}
            where $s$ is a cohesive strength, chosen as $s=1$ for simplicity. Finally, $\bm b_a^\mathrm{solid}$ denotes fluid\nd solid interactions, specified in section~\ref{sec:coupling}.

\section{Coupling scheme}\label{sec:coupling}

    We propose a \emph{virtual interface point} strategy for coupling SPH and LS-DEM, assuming soft contact between the fluid particles and the grains. In this setting, whenever a fluid particle overlaps a grain, LS-DEM contact detection provides the overlapping distance and contact normal, from which a point on the grain surface is identified. This surface point is then treated as a temporary SPH interaction location (not stored as a particle), so that fluid\nd solid pressure-gradient and viscous contributions are evaluated using the same pairwise SPH formulas as for fluid\nd fluid interactions. By action\nd reaction, the resulting interaction is accumulated as force and moment on the grains.

    \subsection{Solid-to-fluid coupling}

\begin{figure}[t!]
        \centering
        \includegraphics[width=0.65\textwidth]{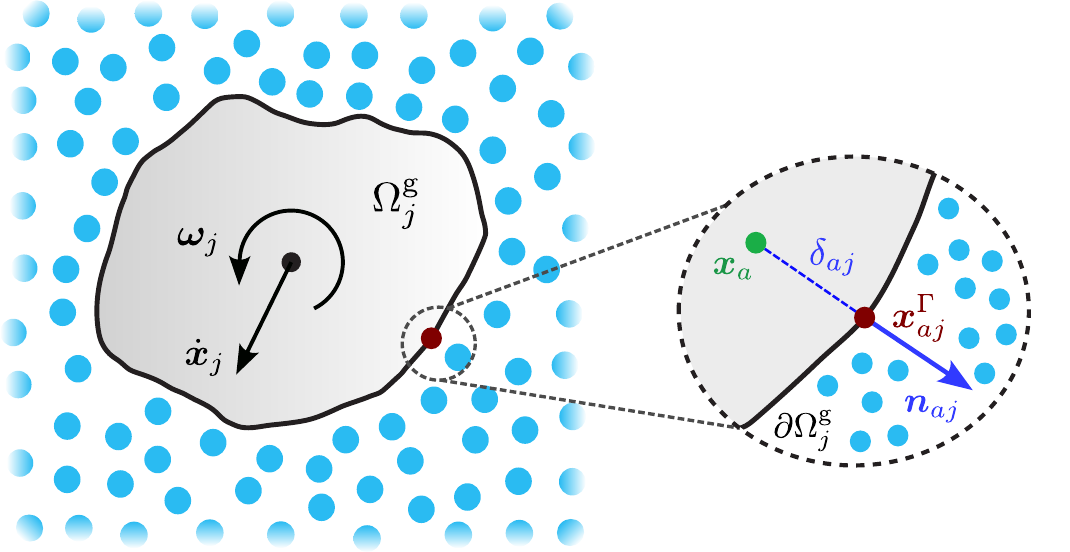}
        \caption{Fluid\nd solid coupling scheme, showing a fluid particle at $\bm x_a$ (green) overlapping grain $j$. The LS-DEM contact detection yields $\delta_{aj}$ and $\bm n_{aj}$, defining a virtual interface point $\bm x_{aj}^{\Gamma}\in\Gamma_j^{\mathrm g}$ (dark red). This point is used to evaluate the fluid\nd solid interaction terms~\eqref{eq:b_aj_solid_def}\nd \eqref{eq:visc_pair_aj}, with grain kinematics given by $\dot{\bm x}_j$ and $\bm\omega_j$.}
        \label{fig:LevelSetVSfluid}
    \end{figure}

    For a fluid particle $a$ interacting with grain $j$, let
    \begin{equation}
        \bm\xi_{aj}
        =
        \mathbf R_j^{\mathsf T}
        \bigl(\bm x_a-\bm x_j\bigr)
    \end{equation}
    denote the particle position in the body-fixed frame of the grain. Standard LS-DEM contact detection provides the penetration distance and spatial outward normal,
    \begin{equation}
        \delta_{aj}
        =
        \max\{-\varphi_j(\bm\xi_{aj}),0\},
        \qquad
        \bm n_{aj}
        =
        \mathbf R_j
        \frac{\nabla\varphi_j(\bm\xi_{aj})}
        {\|\nabla\varphi_j(\bm\xi_{aj})\|}.
    \end{equation}
    The corresponding virtual interface point $\bm x_{aj}^{\Gamma}\in\Gamma_j^{\mathrm g}$ is defined as (Figure~\ref{fig:LevelSetVSfluid})
    \begin{equation}\label{eq:virtualsurfacepoint}
        \bm x_{aj}^{\Gamma}
        \coloneqq
        \bm x_a-\varphi_j(\bm\xi_{aj})\,\bm n_{aj}.
    \end{equation}
    This point provides a discrete location on the fluid\nd solid interface $\Gamma^{\mathrm{fg}}(t)$. We further denote the set of grains interacting with fluid particle $a$ as
    $$
        \mathcal C_a^{\mathrm g}
        \coloneqq
        \{\,j\in\{1,\dots,N^{\mathrm g}\} \ : \ \delta_{aj}>0\,\}.
    $$
    Note that a fluid particle should be contained in at most one grain; however, since grains may overlap slightly in the LS-DEM setting, $\mathcal C_a^{\mathrm g}$ may contain, in principle, more than one index.

    Treating $\bm x_{aj}^{\Gamma}$ as a virtual interface point, we evaluate the fluid\nd solid interaction by inserting the pair $(\bm x_a,\bm x_{aj}^{\Gamma})$ into the SPH expressions~\eqref{eq:pressureAccelerationdiscrete}--\eqref{eq:viscosityAcceleration}. For each contacting pair $(a,j)$, we define the solid-induced acceleration contribution as
    \begin{equation}\label{eq:b_aj_solid_def}
        \bm b_a^{\mathrm{solid}}
        =
        \sum_{j\in\mathcal C_a^{\mathrm g}}
        \bm b_{aj}^{\mathrm{solid}},
        \qquad
        \bm b_{aj}^{\mathrm{solid}}
        =
        -\left(\frac{1}{\rho}\nabla p\right)_{\!aj}
        +
        \nu\,(\nabla^2\bm v)_{aj}.
    \end{equation}
    The pressure-gradient term is evaluated using the symmetric discretization~\eqref{eq:pressureAccelerationdiscrete},
    \begin{equation}\label{eq:press_pair_aj}
        \left(\frac{1}{\rho}\nabla p\right)_{\!aj}
        \approx
        m_{aj}^{\Gamma}
        \left(
            \frac{p_{aj}^{\Gamma}}{\rho_0^2}
            +
            \frac{p_a}{\rho_a^2}
        \right)
        \nabla W(\bm x_a-\bm x_{aj}^{\Gamma},h),
    \end{equation}
    and the viscous term is evaluated using the Brookshaw-type discretization~\eqref{eq:viscosityAcceleration},
    \begin{equation}\label{eq:visc_pair_aj}
        \nu(\nabla^2\bm v)_{aj}
        \approx
        2\nu(n+2)\frac{m_{aj}^{\Gamma}}{\rho_0}
        \left(
            \frac{
                \bigl(\bm v_a-\bm v_j(\bm x_{aj}^{\Gamma})\bigr)
                \cdot
                (\bm x_a-\bm x_{aj}^{\Gamma})
            }{
                \|\bm x_a-\bm x_{aj}^{\Gamma}\|^2+\epsilon h^2
            }
        \right)
        \nabla W(\bm x_a-\bm x_{aj}^{\Gamma},h).
    \end{equation}
    Here, $m_{aj}^{\Gamma}$ is viewed as a numerical interaction weight associated with the virtual interface point and is taken equal to the fluid-particle mass. The virtual pressure is prescribed as
    \begin{equation}
        p_{aj}^{\Gamma}
        =
        k_{\mathrm e}\left|\rho_j-\rho_a\right|,
    \end{equation}
    where $\rho_j$ denotes the density of grain $j$. Thus, $p_{aj}^{\Gamma}$ acts as a repulsive penalty pressure preventing penetration of the fluid particle into the grain. The solid velocity at the virtual interface point follows from rigid-body kinematics,
    \begin{equation}
        \bm v_j(\bm x_{aj}^{\Gamma})
        =
        \dot{\bm x}_j
        +
        \bm\omega_j\times
        (\bm x_{aj}^{\Gamma}-\bm x_j).
    \end{equation}
    The viscous contribution in~\eqref{eq:visc_pair_aj} is driven by the relative velocity at the interface, $\bm v_a-\bm v_j(\bm x_{aj}^{\Gamma})$, which governs the dissipative momentum exchange between the fluid and the moving grain boundary. We note that equation~\eqref{eq:visc_pair_aj} transmits only normal viscous forces. These forces can generate torques and transmit rotational motion to the fluid, particularly for irregular grains. Tangential viscous forces could also be incorporated explicitly by modifying the placement of virtual interaction points or the viscous interaction law, without changing the overall framework.
    
    \subsection{Fluid-to-solid coupling}
    By Newton's third law, the force and moment exerted by the fluid on grain $j$ are obtained from the accumulated action\nd reaction contributions. We define the fluid-contact set of grain $j$ as
    $$
    \mathcal C_j^{\mathrm f}\coloneqq\{\,a\in\{1,\dots,N^{\mathrm f}\} \ : \ \delta_{aj}>0\,\}.
    $$
    The resulting fluid force and moment are, then,
    \begin{equation}\label{eq:fmfluid}
        \bm f_j^{\mathrm{fluid}}
        = -\sum_{a\in\mathcal C_j^{\mathrm f}} m_a\,\bm b_{aj}^{\mathrm{solid}},
        \qquad
        \bm m_j^{\mathrm{fluid}}
        = -\sum_{a\in\mathcal C_j^{\mathrm f}} \bigl(\bm x_{aj}^{\Gamma}-\bm x_j\bigr)\times\bigl(m_a\,\bm b_{aj}^{\mathrm{solid}}\bigr),        
    \end{equation}
    which directly enter the grain equations of motion~\eqref{eq:solid_dynamics}.
    
    \subsection{Coupling summary}

    The proposed fluid\nd solid coupling leverages the Lagrangian setting of both LS-DEM and SPH, and employs the contact detection of LS-DEM to impose moving solid boundaries within the WCSPH formulation. The coupled system is then advanced explicitly in time using a common step $\Delta t$, resulting in a convenient implementation scheme. We summarize the overall procedure as follows. At each time step:

    \begin{enumerate}[itemsep=0.25em]
        \item Compute the WCSPH variables from~\eqref{eq:discreteDensitya}\nd \eqref{eq:pressurea} and assemble the fluid\nd fluid contributions in the discrete momentum balance~\eqref{eq:momentumFluid} using~\eqref{eq:pressureAccelerationdiscrete} and~\eqref{eq:viscosityAcceleration}.
        \item Use LS-DEM contact detection to identify overlapping fluid\nd solid pairs, defining the contact sets $\mathcal C_a^{\mathrm g}$ and $\mathcal C_j^{\mathrm f}$ and constructing the corresponding virtual interface points~\eqref{eq:virtualsurfacepoint}.
        \item Evaluate solid-to-fluid coupling by computing $\bm b_{aj}^{\mathrm{solid}}$ and $\bm b_a^{\mathrm{solid}}$ via~\eqref{eq:b_aj_solid_def}\nd \eqref{eq:visc_pair_aj}. Accumulate the corresponding fluid-to-solid force and moment on each grain via~\eqref{eq:fmfluid}.
        \item Compute grain--grain and grain--wall contact forces and moments using~\eqref{eq:contactForce} and~\eqref{eq:momentni}.
        \item Advance SPH and LS-DEM states with time step $\Delta t$ by integrating~\eqref{eq:momentumFluid} and~\eqref{eq:solid_dynamics}\nd \eqref{eq:quatkin}.
    \end{enumerate}

\section{Three-dimensional simulations}\label{sec:numerical}

    In this section, the performance and robustness of the proposed SPH\nd LS-DEM coupling scheme are evaluated through four numerical benchmarks. The tests are designed to assess the ability of the method to handle fluid\nd solid interactions under complex boundary conditions and different physical scenarios.

    \subsection{Dam-break flow impacting an obstacle} \label{subsec:dambreak}
    
    This first test consists of a classic dam-break problem, which is a standard benchmark for validating free-surface dynamics. Specifically, we assess numerical stability and the ability of the proposed method to capture the dynamics of fluid fronts after impacting an obstacle and a wall.

      \begin{figure}[t!]
        \centering
        \small
        \includeinkscape[width=0.8\textwidth]{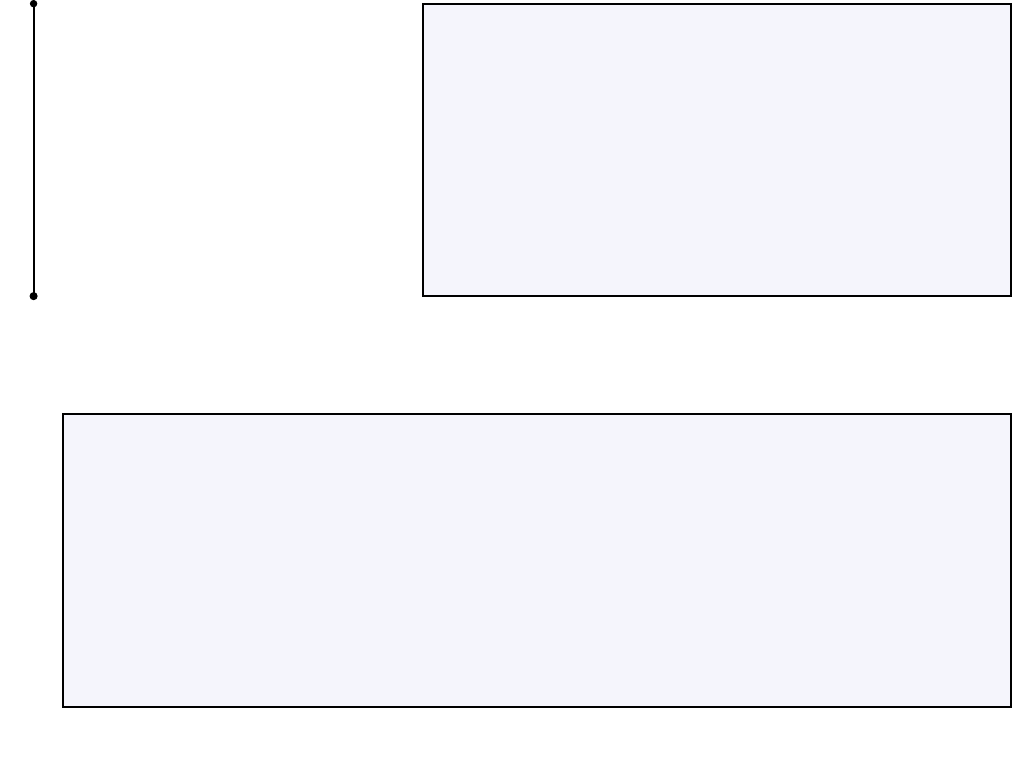}     
        \caption{Schematic representation of initial conditions for the dam-break test: plan view (top) and elevation view (bottom) of the tank. The fluid volume (light blue) is initially contained by a rigid wall and then subsequently released, allowing it to flow and impact the obstacle (dark gray) and the rightmost wall. Dimensions given in m.}
        \label{fig:damBreakTestDims}
    \end{figure} 

    \begin{table}[b!]
        \centering
        \caption{Simulation parameters for the dam-break problem.}        
        \begin{tabular}{l l l }
        \toprule
        \textbf{Simulation parameter} & \textbf{Value}  & \textbf{Units}          
        \\ 
        \midrule
        Simulation time  & 6 & s
        \\
        Stiffness constant (SPH) $k_\mathrm{e}$ & 40 & m$^{2}$/s$^{2}$
        \\
        Fluid particle mass $m$ & 0.02 & kg 
        \\ 
        Smoothing length $h$ & 0.1 & m
        \\
        Time increment $\Delta t$  & 0.0001 & s
        \\ 
        Target density $\rho_{0}$ & 1000 & kg/m$^{3}$
        \\ 
        Gravity $g$  & 9.81 & m/s$^{2}$
        \\ 
        Kinematic viscosity $\nu$  & 0.000001 & m$^{2}$/s
        \\ 
        \bottomrule      
        \end{tabular}
        \label{table:dambreakSimmulationParameters}
    \end{table}

    \begin{figure}[t!]
        \centering
        \small
        \includegraphics[width=0.9\textwidth]{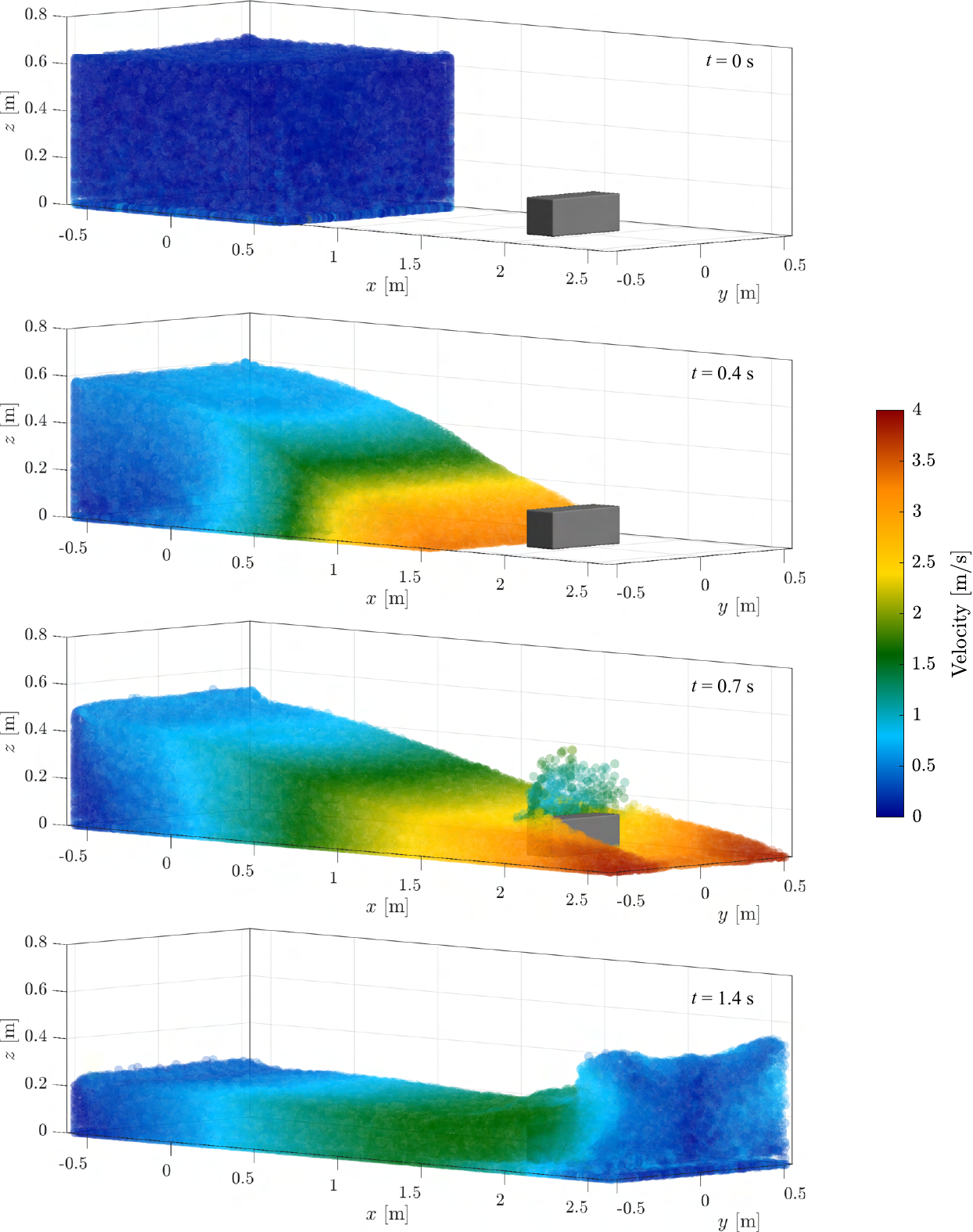}
        \caption{Snapshots of the dam-break flow impacting an obstacle.}
        \label{fig:bvp1}
    \end{figure} 
    
    \begin{figure}[t!]
        \centering
        \includegraphics[width=0.95\textwidth]{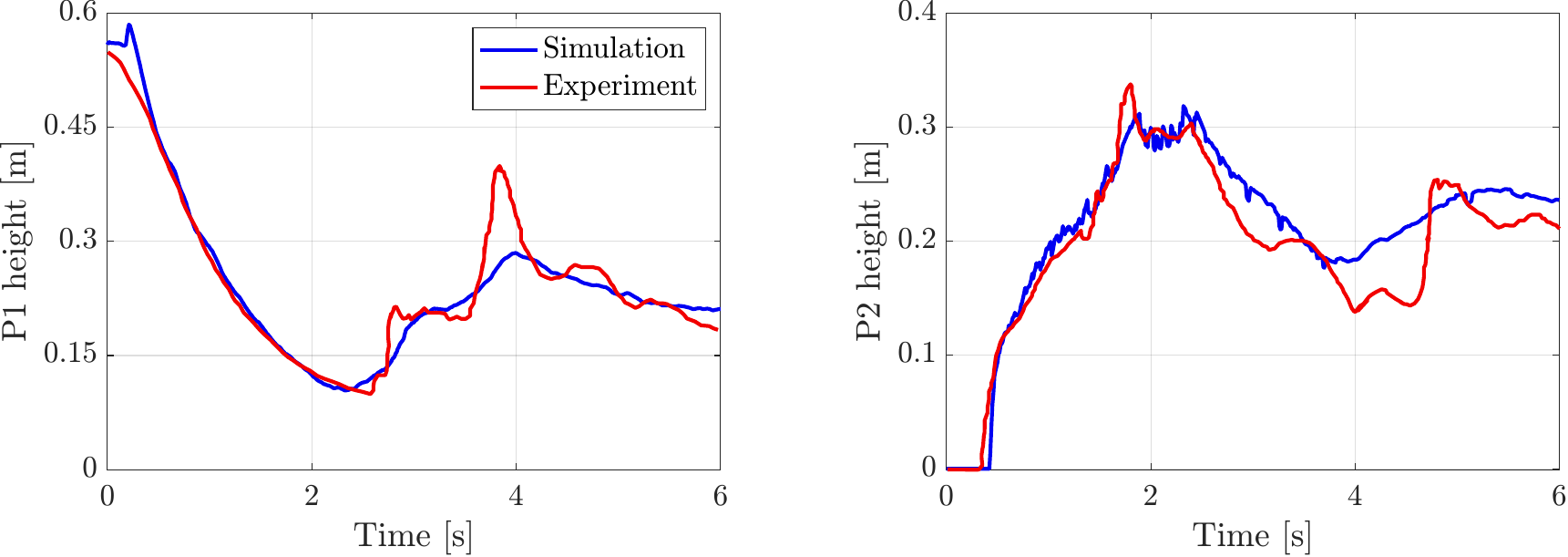}
        \begin{picture}(0,0)
        \put(-442,5){(a)}
        \put(-205,5){(b)}
        \end{picture}       
        \caption{Fluid heights at probes P1 and P2, comparing the SPH\nd LS-DEM simulation with experimental data~\cite{kleefsman2005volume}.}
        \label{fig:H2H4}
    \end{figure} 

    The setup is based on an experiment conducted at the Maritime Research Institute Netherlands (MARIN)~\cite{kleefsman2005volume}. An open-top tank with dimensions 3.22$\times$1$\times$1~m is divided into two sections, with the left section initially filled with water to a depth of 0.55~m and separated from the remainder of the tank by a vertical wall. An obstacle consisting of a rigid block is placed in the right section and fixed to the ground. The water is released by rapidly removing the wall, triggering flow towards the obstacle and the opposite wall. Height probes measure the water level, with probes P1 and P2 (corresponding to H4 and H2 in the experiment) located in the reservoir and in front of the obstacle. Figure~\ref{fig:damBreakTestDims} shows the geometry and dimensions of the system, as well as the locations of the probes and the~obstacle.

    We conduct a numerical simulation using 27864 SPH fluid particles, which are poured into a reservoir enclosed by rigid walls according to the initial conditions in Figure~\ref{fig:damBreakTestDims}.  The obstacle is modeled as an LS-DEM solid with a prismatic shape. The fluid is allowed to settle and converge to a state of near-zero kinetic energy. The spiky kernel~\eqref{eq:spikyK} is employed, and the simulation parameters presented in Table~\ref{table:dambreakSimmulationParameters} are chosen to approximate the behavior of water. In particular, the artificial stiffness $k_\mathrm{e}$ inherent to WCSPH is chosen to ensure a suitable trade-off between accuracy and numerical stability, with too-low values deviating strongly from the incompressibility condition and too-high values leading to numerical~instability.

    Figure~\ref{fig:bvp1} shows snapshots of the fluid flow over a period of 1.4 s. The fluid is at rest at $t=0$, and the simulation time starts when the wall is removed. Then, the SPH particles are released and driven by gravity. After approximately $0.4$ s, the fluid impacts the obstacle and splashes upward, forming a wave crest. The flow continues until the fluid hits the back wall, generating a backward impact wave.

    Figure~\ref{fig:H2H4} shows the fluid heights at the probe locations P1 and P2 and compares the numerical simulation with the experimental benchmark~\cite{kleefsman2005volume}. In the simulation, the height at each probe is computed from the positions of the SPH particles. To filter out free-surface noise, we average the vertical coordinates of the particles in the uppermost layer within a narrow sampling radius centered at the probe. The comparisons reveal excellent agreement between the experimental results and the proposed SPH\nd LS-DEM method across all stages. Notably, the agreement persists as the fluid front crosses the obstacle, hits the back wall, travels back past the probes (when the fluid height at P2 is largest), returns to the reservoir, and finally advances again, passing probe P2 for a third time at about $5$~s.

    \subsection{Darcy-type flow through granular packings}
    \label{subsec:kozeny}

    The second example addresses fluid flow through a porous granular medium, considering packings of spheres and cubes. The response predicted by the present method under nearly steady saturated flow is compared with Darcy-type behavior and the Kozeny\nd Carman equation, a classical empirical relation used to evaluate the pressure drop $\Delta p$ as fluid flows through a granular packing. Specifically,
    \begin{equation}
        \frac{\Delta p}{L} = \frac{150\,\eta}{\Phi_{\rm s}^2 \, d_{\rm p}^2} \frac{(1 - \varepsilon)^2}{\varepsilon^3} V_0,
        \label{eq:KC}
    \end{equation}
    where $L$ is the packing length, $\varepsilon$ is the packing porosity, $d_{\rm p}$ is the diameter of a sphere with volume equal to the typical grain in the packing, $\Phi_{\rm s}$ is the sphericity of the typical grain in the packing, and $V_0$ is the fluid's superficial velocity. This expression represents a Darcy-type law, where the effective permeability is expressed in terms of the physical properties of the granular medium.

        \begin{table}[b!]
        \centering
        \caption{Simulation parameters for the fluid phase (sections \ref{subsec:kozeny}, \ref{subsec:granularCollapse}, and \ref{subsec:cyclicLoading}).}   
        \begin{tabular}{l l l}
        \toprule
        \textbf{Simulation parameter} & \textbf{Value}   & \textbf{Units}       
        \\ 
        \midrule
        Simulation time  & 6 & s
        \\
        Fluid particle mass $m$ & 0.01 & kg 
        \\ 
        Smoothing length $h$ & 0.05 & m
        \\
        Time increment $\Delta t$  & 0.00001 & s
        \\ 
        Target density $\rho_{0}$ & 1000 & kg/m$^{3}$
        \\ 
        Gravity $g$  & 9.81 & m/s$^{2}$
        \\ 
        Kinematic viscosity $\nu$  & 0.00005 & m$^{2}$/ s
        \\
        \bottomrule      
        \end{tabular}
        \label{table:flowParams}
    \end{table}

    \begin{figure}[t!]
        \centering
        \small
        \includegraphics[width=0.46\textwidth]{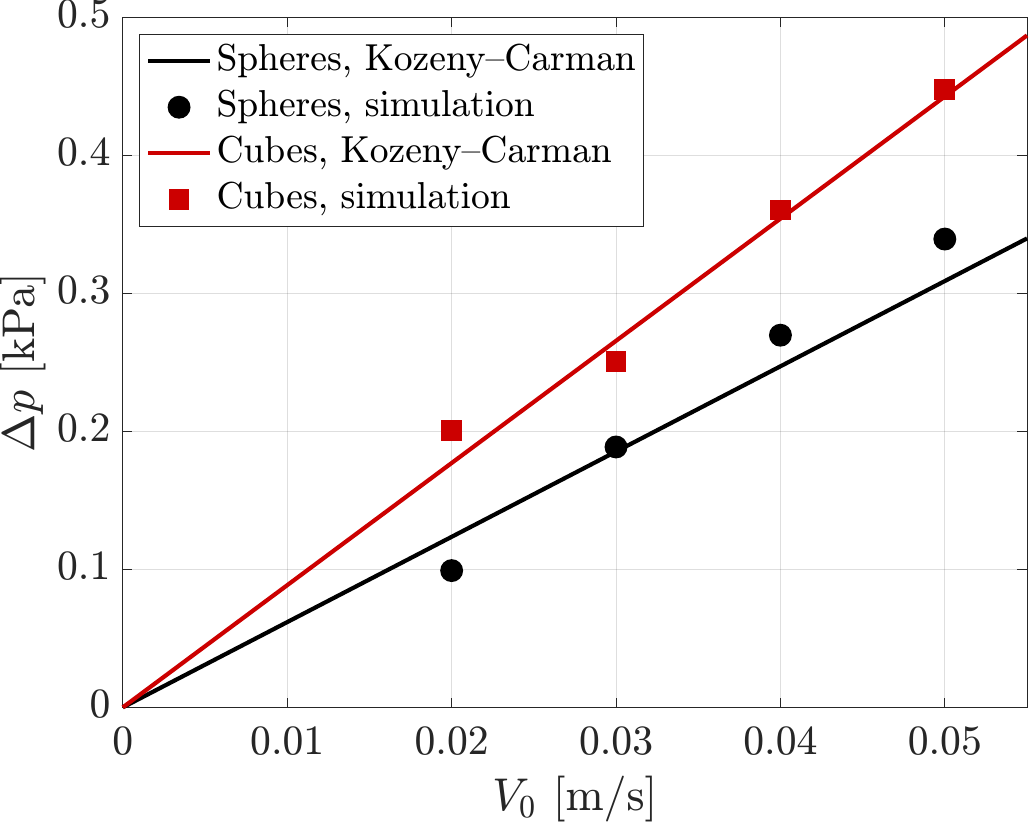} 
        \caption{Pressure drop versus fluid velocity through packings of spheres and cubes, comparing the SPH\nd LS-DEM simulations with the Kozeny\nd Carman relation~\eqref{eq:KC}.}
        \label{fig:Kozeny}
    \end{figure}

    \begin{figure}[t!]
        \centering
        \includegraphics[width=0.95\textwidth]{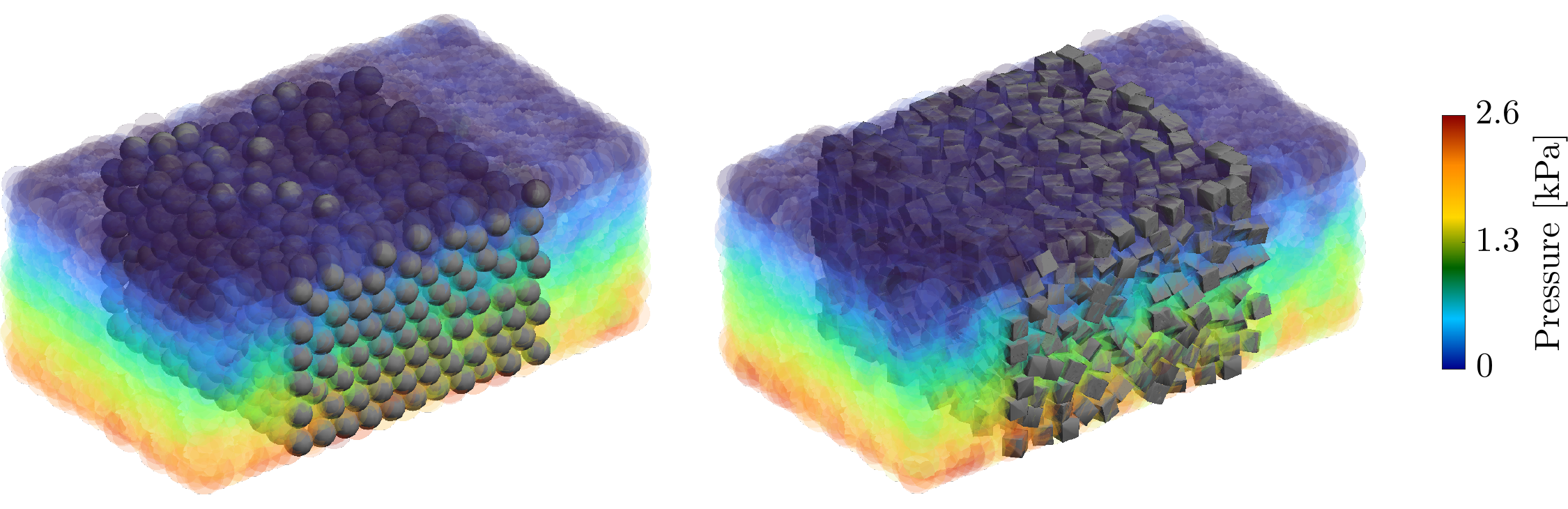}
        \begin{picture}(0,0)
        \put(-438,9){(a)}
        \put(-238,9){(b)}        
        \end{picture}  
        \caption{Snapshots of the pressure field during flow through (a) spheres and (b) cubes in the nearly stationary~regime.}
        \label{fig:kc_snaps}
    \end{figure}     

    The computational domain consists of a $2.1$ m horizontal channel with a cross-sectional area $0.2412$~m$^2$, where the granular packing occupies the region $0.8\leq x \leq 1.4 $. The cubic spline kernel~\eqref{eq:cubic} is employed, together with the simulation parameters in Table~\ref{table:flowParams}. In addition, $k_{\rm e}$ is set to 40~m$^{2}$/s$^{2}$. The channel, including the pore space, is initially filled with $55872$ SPH particles, and the fluid is allowed to reach hydrostatic equilibrium under gravity. A piston placed at $x=0$ is then moved in the positive $x$-direction at a constant velocity $V_0$, while the outlet remains open. The pressure drop across the packing is evaluated as $\Delta p=\bar{p}_{\mathrm{up}}-\bar{p}_{\mathrm{down}}$, where $\bar{p}_{\mathrm{up}}$ and $\bar{p}_{\mathrm{down}}$ denote cross-sectional pressure averages upstream and downstream of the packing, respectively. Because both sections span the same vertical range, the hydrostatic contribution cancels. The fluid velocity is similarly averaged over the cross-section at several locations within the packing. Once the pressure drop and velocity averages reach an approximately stationary state, time averages are computed, while verifying that the cross-sectional velocity averages remain close to $V_0$.
    
    Figure~\ref{fig:Kozeny} compares the simulated pressure drop with the Kozeny\nd Carman equation~\eqref{eq:KC} for packings of spheres ($\varepsilon = 0.4319$, $d_{\rm p} = 0.054$~m, $\Phi_{\rm s} = 1$) and cubes ($\varepsilon = 0.4402$, $d_{\rm p} = 0.054$~m, $\Phi_{\rm s} = 0.8$), while Figure~\ref{fig:kc_snaps} shows representative pressure fields during the nearly stationary regime. For both spheres and cubes, $\Delta p$ increases approximately linearly with $V_0$, consistent with a Darcy-type response. Moreover, the packing of cubes produces a larger pressure drop than the packing of spheres at the same velocity, indicating a lower effective permeability. The simulated values closely follow equation~\eqref{eq:KC},  showing that the SPH\nd LS-DEM framework reproduces the dependence of pressure drop on fluid velocity and grain shape, in agreement with the Kozeny\nd Carman~equation.

    \subsection{Granular column collapse}
    \label{subsec:granularCollapse}

\begin{figure}[t!]
        \centering
        \small
        \includegraphics[width=0.95\textwidth]{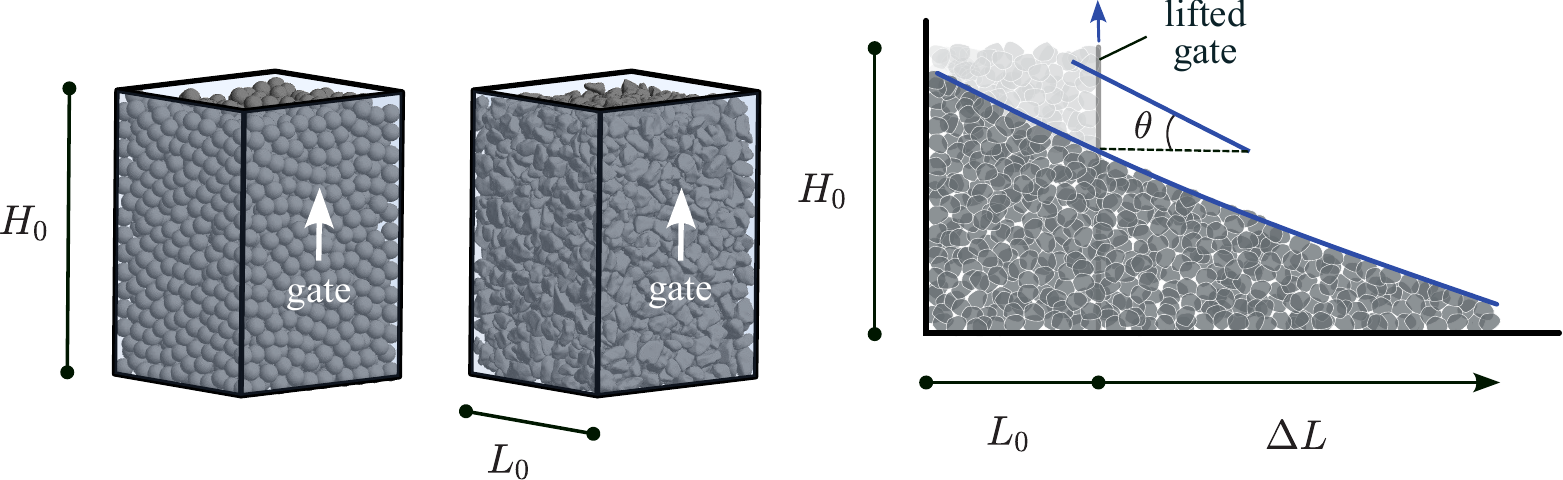}
        \begin{picture}(0,0)
        \put(-440,5){(a)}
        \put(-212,5){(b)}        
        \end{picture}  
        \caption{Schematic representation of (a) the 3D granular column setup with initial dimensions of the column $H_{0}=1$~m and $L_{0}=0.75$~m for both spheres and angular grains and (b) the collapse process with runout distance $\Delta L$.}
        \label{fig:collapse_scheme}
    \end{figure}

    The third example considers a series of three-dimensional granular-column collapse tests involving dry and initially saturated assemblies (Figure~\ref{fig:collapse_scheme}). The tests are carried out separately for assemblies of $2000$ spheres and $1750$ irregular grains corresponding to a regolith simulant~\cite{Kawamoto2016Jun}.
    
    The diameter of the spheres is chosen as $0.07$~m to achieve a packing fraction similar to that of the irregular grains. Table~\ref{table:grainParams} shows the physical properties and simulation parameters for the solid phase. For the fluid phase, the cubic spline kernel~\eqref{eq:cubic} is employed, together with the parameters in Table~\ref{table:flowParams} and $k_{\rm e}= 40$~m$^{2}$/s$^{2}$. In all cases, the column has an initial height $H_0 = 1$~m and base length $L_0 = 0.75$~m, and is released by instantaneously lifting a vertical gate. In the wet configuration, the granular assembly is initially filled with water while the surrounding environment is dry. Thus, unlike a fully submerged collapse, where the surrounding fluid contributes additional buoyancy and drag, fluid effects here are expected to arise primarily from pore pressure changes and drag within the granular assembly~\cite{HeShiYu2022}.

  \begin{table}[t!]
        \centering
        \caption{Simulation parameters for the solid phase (sections \ref{subsec:kozeny}, \ref{subsec:granularCollapse}, and \ref{subsec:cyclicLoading}).}        
        \begin{tabular}{l l l}
        \toprule
        \textbf{Simulation parameter} & \textbf{Value}   & \textbf{Units}       
        \\ 
        \midrule
        Simulation time  & 6 & s
        \\ 
        Normal stiffness constant $k_{\rm n}$ & 1000000 & N/m
        \\ 
        Tangential stiffness constant $k_{\rm s}$ & 900000 & N/m
        \\  
        Time increment $\Delta t$  & 0.00001 & s
        \\ 
        Grain density & 2500 & kg/m$^3$
        \\ 
        Gravity $g$  & 9.81 & m/s$^2$
        \\ 
        Grain\nd grain friction coefficient $\mu$ & 0.65 & --
        \\ 
        \bottomrule
        \end{tabular}
        \vspace*{2em}
        \label{table:grainParams}
    \end{table}

\begin{figure}[b!]
        \centering
        \includegraphics[width=0.9\textwidth]{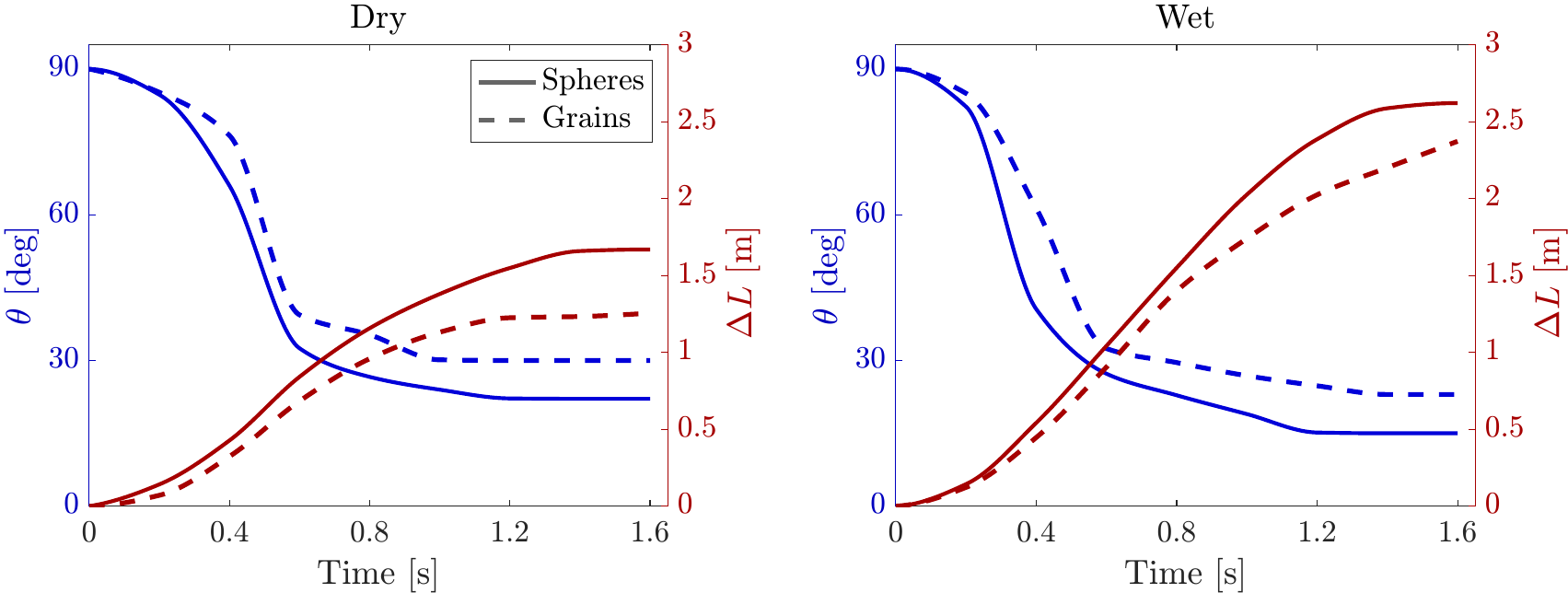}
        \begin{picture}(0,0)
        \put(-442,5){(a)}
        \put(-211,5){(b)}        
        \end{picture}  
        \caption{Time histories of the tangent angle at the opening gate location $\theta$ and runout distance $\Delta L$ (Figure~\ref{fig:collapse_scheme}) for spheres and angular grains in (a) dry conditions and (b) wet conditions (see Figure~\ref{fig:bvp2_dry_wet}).}
        \label{fig:fluidization_curves}
    \end{figure}

\begin{figure}[b!]
        \centering
        \small
        \includegraphics[width=0.95\textwidth]{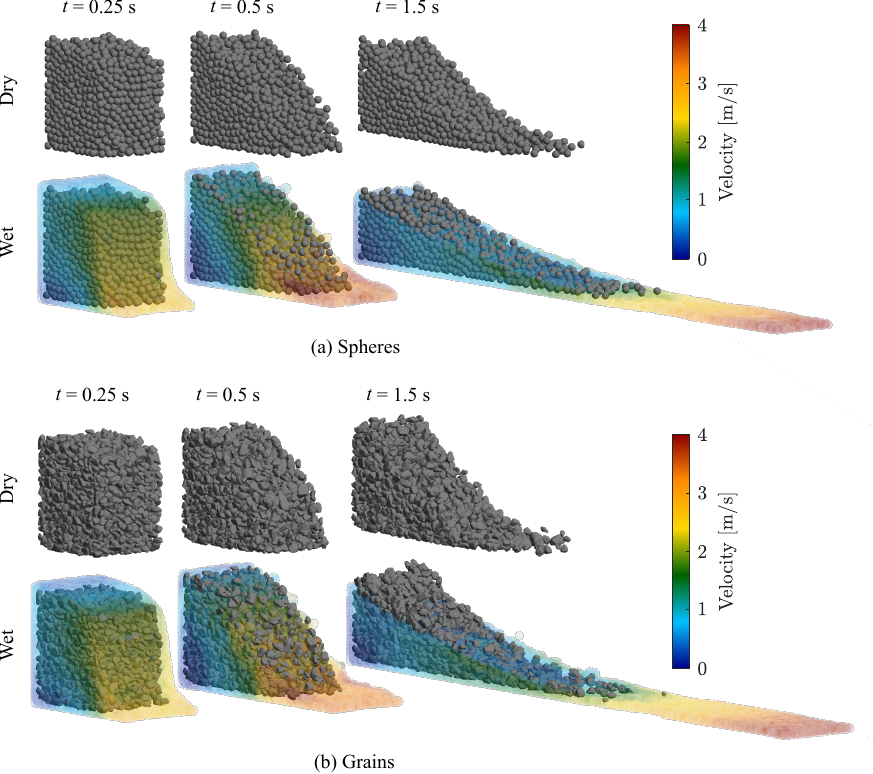}
        \caption{Snapshots of the granular collapse process for (a) spheres and (b) regolith grains in dry and wet conditions.}
        \label{fig:bvp2_dry_wet}
    \end{figure}

    Figure~\ref{fig:fluidization_curves}(a) shows that the free-surface angle at the gate, $\theta(t)$, decreases rapidly after the gate is lifted and as the column relaxes, while the runout length, $\Delta L(t)$, grows monotonically. Consistent with these trends, the snapshots in Figure~\ref{fig:bvp2_dry_wet} show the material spreading horizontally and eventually stabilizing. For the spheres, the front advances faster and the final free surface is noticeably flatter, whereas the irregular grains retain a steeper surface profile. This behavior reflects the reduced geometric interlocking and lower rotational resistance of the spheres, and is in line with previous numerical studies showing that runout increases with particle regularity~\cite{ZhangTahmasebi2022}.

    \begin{table}[t!]
        \centering
        \caption{Comparison of the angle of repose between the simulations and references under wet and dry conditions, for spherical and irregular regolith grains. The reference value for regolith grains is computed according to~\citet[equation~(4)]{cho2006}, using the average roundness of the sampled grains (0.6229).}
        \begin{tabular}{l c c}
        \toprule
        \textbf{Condition/particle type} & \textbf{Simulation [$^\circ$]} & \textbf{Reference [$^\circ$]} \\ 
        \midrule
        Dry spheres & 22.1 & 21 \cite{cho2006} \\ 
        Wet spheres & 15  & 15 \cite{ceccato2020}\\ 
        Dry  regolith grains  & 30 & 31.46 \cite{cho2006} \\ 
        Wet regolith grains  & 23 & -- \\ 
        \bottomrule
        \end{tabular}
        \label{table:angleComparison}
    \end{table}    

    Compared with the dry configuration, Figure~\ref{fig:fluidization_curves}(b) shows that the wet configuration exhibits a lower stabilized angle $\theta$ and a larger runout $\Delta L$. This behavior, associated with fluid\nd grain interactions and reduced effective frictional resistance within the granular skeleton, agrees with the partially immersed column experiments and simulations of~\citet{HeShiYu2022}. Accordingly, the snapshots in Figure~\ref{fig:bvp2_dry_wet} show lower free-surface slopes in the wet case. Furthermore, the contrast between spheres and angular grains is slightly reduced in the presence of water, suggesting that fluid\nd grain interactions mitigate the mobility reduction associated with particle irregularity~\cite{ZhangTahmasebi2022}. The final angles of repose are summarized in Table~\ref{table:angleComparison}, showing close agreement with the available reference values for both particle shapes and fluid conditions.
    
    \subsection{Cyclic triaxial loading of a granular specimen} \label{subsec:cyclicLoading}

    The last example addresses the cyclic behavior of granular soils, a canonical problem in geomechanics with direct relevance to liquefaction. We illustrate the ability of the proposed method to capture liquefaction-like behavior through triaxial cyclic loading simulations under both dry and saturated undrained conditions. The objective is to capture the expected difference in mechanical behavior between dry and fully saturated cases, focusing on the cumulative effects of cyclic loading on effective stress paths. 
    
\begin{figure}[b!]
        \centering
        \includegraphics[width=0.9\textwidth]{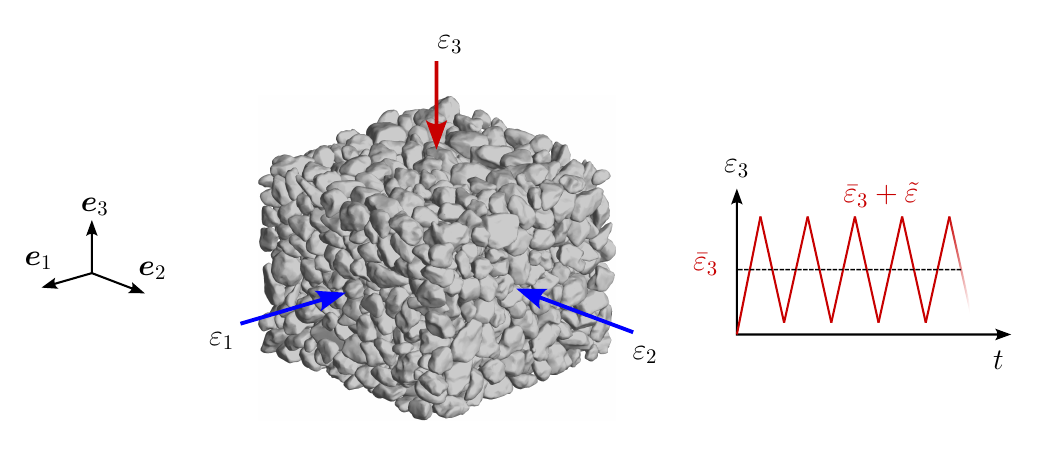}  
        \caption{Granular specimen and loading conditions for the cyclic loading test.}
        \label{fig:bvp3}
    \end{figure}

    The specimen was prepared by pluviating $1000$ irregularly shaped regolith grains~\cite{Kawamoto2016Jun} into a rigid-walled container, resulting in an assembly of dimensions $0.8 \times 0.8 \times 0.65$ m, with the shortest dimension in the $\bm e_3$ direction (Figure~\ref{fig:bvp3}). The specimen was then isotropically compressed to a confining pressure of $6.9$ kPa. For the saturated case, $34438$ SPH fluid particles were allowed to flow into the pores, and the coupled system was allowed to reach static equilibrium at the same effective stress as the dry sample. Constant-amplitude axial-strain cycles were imposed in both cases, with $\bar{\varepsilon}_3=0.01 $ and $\tilde{\varepsilon}_3\in [-0.01,0.01]$. Meanwhile, the lateral confining stress was kept constant following a servo-controlled scheme~\cite{o2011particulate}. Global damping of $100$~s$^{-1}$ was applied to maintain quasi-static conditions, with an inertial number below $10^{-3}$~\cite{RadjaiDubois2011}. As in the previous examples, we use the parameters in Table~\ref{table:grainParams} for the grains and Table \ref{table:flowParams} for the fluid, together with the cubic spline kernel~\eqref{eq:cubic} and $k_{\rm e}=400$~m$^2$/s$^2$ .



    Pore pressure can significantly affect the mechanical response of granular packings, even at low strain rates. During undrained shearing, a contractive tendency of the packing can generate positive excess pore pressure, reducing the effective stress carried by grain contacts and, consequently, their frictional resistance. Following Terzaghi's effective stress principle, the average total stress of the sample can be decomposed as
    \begin{equation}
        \bar{\bm{\sigma}} = \bm{\sigma}' + P \mathbf{I},
    \end{equation}
    where $P$ is the average pore pressure in the fluid and $\bm{\sigma}'$ is the effective stress in the granular phase, computed from the contact forces and branch vectors as (cf.~\citet{andrade2022})
    \begin{equation}
        \bm{\sigma}'
        = \frac{1}{V}\sum_{c\in\mathcal{C}}
        \operatorname{sym} \left(\bm{f}^{c}\otimes\bm{\ell}^{c}\right),
    \end{equation}
    where $V$ is the specimen volume, $\mathcal{C}$ is the set of contacts spanning the entire assembly (counted once), $\bm{f}^{c}$ is the contact force, and $\bm{\ell}^{c}$ is the corresponding branch vector. The orientations of $\bm{f}^{c}$ and $\bm{\ell}^{c}$ are chosen consistently with the compression-positive stress convention in the present example. The effective stress $\bm{\sigma}'$ is conveniently characterized by the mean effective stress $p'$ and the equivalent deviatoric stress $q$, defined as
    \begin{equation}
        p'=\displaystyle \frac{\mathrm{tr}(\bm{\sigma}')}{3}, \qquad q=\displaystyle \sqrt{\frac{3}{2}(\operatorname{dev}\bm{\sigma}':\operatorname{dev}\bm{\sigma}')}.
    \end{equation}

\begin{figure}[t!]
        \centering
        \includegraphics[width=0.9\textwidth]{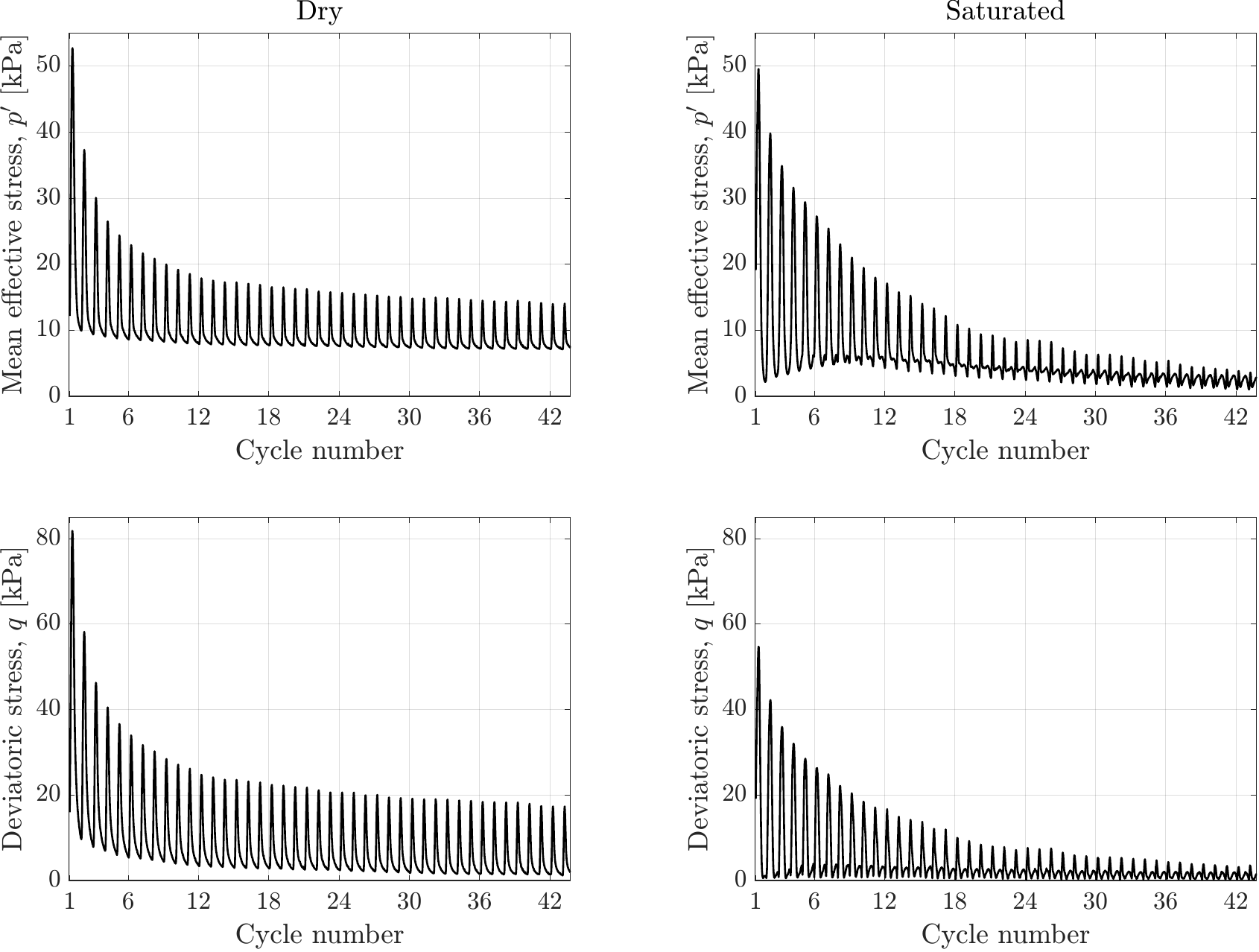}
        \begin{picture}(0,0)
        \put(-428,164){(a)}
        \put(-205,164){(b)}
        \put(-428,5){(c)}
        \put(-205,5){(d)}        
        \end{picture}  
        \caption{Evolution of mean effective stress $p'$ and deviatoric stress $q$ for the dry and saturated cases.}
        \label{fig:mean_dev_stress}
    \end{figure}

    Before cyclic loading, both the dry and saturated specimens are in isotropic equilibrium at the same initial mean effective stress. As shown in Figure~\ref{fig:mean_dev_stress}(a,c), the dry specimen exhibits an initial reduction in mean effective stress $p'$ and deviatoric stress amplitude, followed by approximately stable cyclic oscillations. By contrast, as illustrated in~Figure~\ref{fig:mean_dev_stress}(b,d), the saturated specimen exhibits a cyclic liquefaction-like response, characterized by a progressive reduction in mean effective stress and deviatoric stress.
    
       \begin{figure}[t!]
        \centering
        \includegraphics[width=0.9\textwidth]{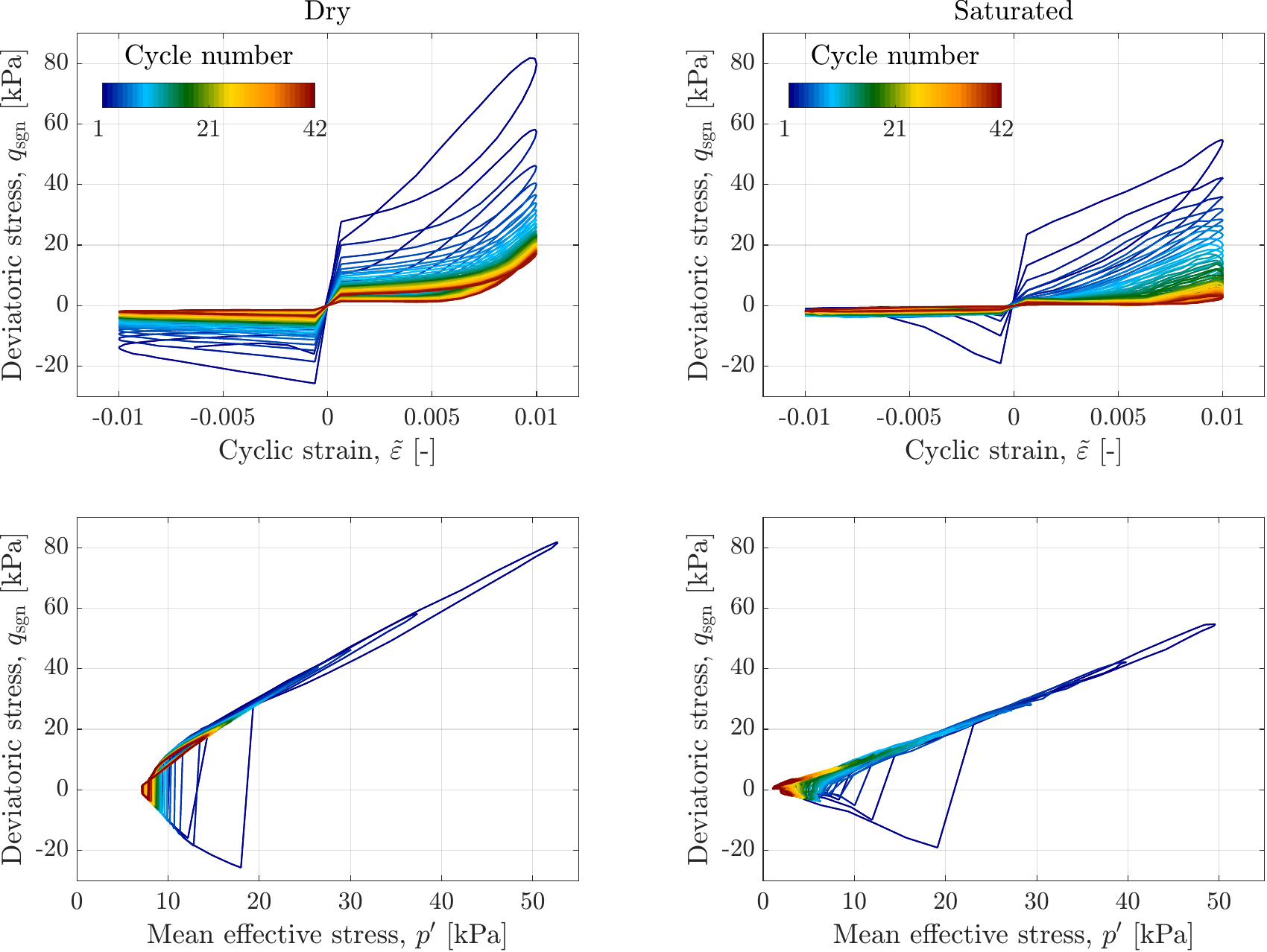}
        \begin{picture}(0,0)
        \put(-428,164){(a)}
        \put(-205,164){(b)}
        \put(-428,5){(c)}
        \put(-205,5){(d)}         
        \end{picture}  
        \caption{Deviatoric stress response, $q_{\mathrm{sgn}}\coloneqq q\,\operatorname{sign}(\tilde{\varepsilon}_3)$, as a function of (a,b) cyclic strain and (c,d) mean effective~stress.}
        \label{fig:phase_space}
    \end{figure}          
    
    Figure~\ref{fig:phase_space} shows the deviatoric stress response as a function of cyclic strain (a,b) and mean effective stress (c,d), where $q_{\mathrm{sgn}}\coloneqq q\,\operatorname{sign}(\tilde{\varepsilon}_3)$ distinguishes the positive and negative halves of the strain cycle. During the first few cycles, both dry and saturated specimens exhibit relatively high stiffness. With continued cycling, the loops progressively flatten in both cases (cyclic softening), but the dry specimen stabilizes while the saturated specimen exhibits a continued reduction in deviatoric stress amplitude. The method thus captures the progressive build-up of excess pore pressure under undrained cyclic loading and the associated loss of effective stress and shear resistance.
    
\section{Conclusions and perspectives}\label{sec:conc}

    We have presented a coupled SPH\nd LS-DEM method for fluid\nd solid interaction. The method naturally handles geometrically complex, moving interfaces, as encountered in granular systems under saturated or partially saturated conditions. The framework is particularly advantageous due to the natural compatibility of the Lagrangian descriptions adopted for the fluid and solid phases. Moreover, the proposed formulation and computational implementation enable a simple fluid\nd solid interaction procedure by leveraging the level-set function discretization inherent to LS-DEM. 
    
    Several examples were presented to highlight the ability of the framework to capture complex behaviors across a range of physical systems. A dam-break problem with flow impacting an obstacle was first considered and compared against experimental observations, showing close agreement. A Darcy-type flow through packings of spheres and cubes was then examined, demonstrating close agreement with the Kozeny\nd Carman relation. The third example addressed a granular collapse test to compare the dynamics of dry assemblies and partially saturated configurations, as well as the response of spherical grains versus irregular particles from a regolith simulant. Finally, the last example demonstrated an application to cyclic loading of granular media, in which a fluid\nd granular assembly qualitatively reproduced key features of liquefaction-like phenomena.

    We close by discussing open issues of the proposed model, to be addressed in future work:

    \begin{itemize}
        \item The choice of SPH kernel affects the accuracy and stability of the function approximations. In this work, the spiky and cubic spline kernels are used; however, alternative kernels that provide more accurate approximations of the spatial derivatives could be explored.
        \item In this work, we use a WCSPH scheme to estimate the pressure from density deviations, which entails treating the fluid as weakly compressible. A closer approximation of incompressibility could be obtained by solving a Poisson equation for the pressure.
        \item The viscous forces are evaluated explicitly using the approximation proposed by~\citet{monaghan2005sph}, which is stable for low viscosities. To simulate highly viscous fluids, implicit methods can instead be employed.
        \item Finally, computational cost could be reduced through GPU acceleration, for example using CUDA, enabling finer spatial discretizations of the fluid and larger granular assemblies.
    \end{itemize}     

    \noindent Furthermore, the proposed framework sets the stage for various multiphysics systems, in which the SPH particles and the contact evaluation procedure can be used to describe different physical interactions.

\section*{Acknowledgements}

\noindent This work was supported by the U.S.~Department of Energy's Advanced Simulation and Computing (ASC) Physics and Engineering Models (PEM) Enabling Manufacturing (EM) through Los Alamos National Laboratory and LANL--Michigan SPARC. Los Alamos National Laboratory is operated by Triad National Security, LLC, for the National Nuclear Security Administration under Contract No.~89233218CNA000001.

\small
\bibliography{references}

\end{document}